\documentclass[12pt,reqno]{amsart}
\usepackage[english]{babel}
\usepackage[margin=1.5in]{geometry}
\usepackage[applemac]{inputenc}
\usepackage{tikz}
\usetikzlibrary{arrows,automata}
\usepackage{times}
\usepackage{tikz-cd}

\newcommand{\C}{{\mathbb C}}

\newcommand{\G}{{\mathcal G}}

\newcommand{\lb}{\lambda}
\newcommand{\td}{{\mathbb T}}

\newcommand{{\pf}}{{\bf Proof. }}

\newtheorem{thm}{Theorem}[section]
\newtheorem{corr}[thm]{Corollary}
\newtheorem{lem}[thm]{Lemma}
\newtheorem{prop}[thm]{Proposition}
\newtheorem{exam}[thm]{Example}

\newcommand{\zb}{$\hfill\Box$}

\makeatletter\@addtoreset{equation}{section} \makeatother

\usetikzlibrary{arrows}

\tikzset{
  treenode/.style = {align=center, inner sep=1pt, text centered,
    font=\sffamily},
  arn_n/.style = {treenode, circle, white, font=\sffamily\bfseries, draw=black,
    fill=black, text width=1.5em},
  arn_r/.style = {treenode, circle, black, draw=black,
    text width=1.5em, very thick},
  arn_x/.style = {treenode, rectangle, draw=black,
    minimum width=0.5em, minimum height=0.5em}
}

\begin{document}

\title{Joint Spectrum and the Infinite Dihedral Group}

\author[R. Grigorchuk]{Rostislav Grigorchuk$^1$}
\footnotemark[1]
\address{Rostislav Grigorchuk: Department of Mathematics, Texas A\& M University,
College Station, TX 77843.} \email{grigorch@math.tamu.edu}

\author[R. Yang]{Rongwei Yang}
\address{Rongwei Yang: Department of Mathematics and Statistics, SUNY at Albany,
Albany, NY 12222, U.S.A.} \email{ryang@math.albany.edu}

\maketitle

\footnotetext[1]{$^1$ The first author is supported by NSA grant H98230-15-1, the
Swiss National Science Foundation, and  ERC AG COMPASP.}

\footnotetext{
\hangindent=1.2em\emph{Key words and phrases}: projective joint spectrum, dihedral group, $C^*$-algebra, weak containment, Maurer-Cartan form, Fuglede-Kadison determinant, group of intermediate growth, self-similar representation, support of representation, dynamics.}

\begin{abstract}
For a tuple $A=(A_1,\ A_2,\ ...,\ A_n)$ of elements in a unital Banach algebra ${\mathcal B}$, its {\em projective joint spectrum} $P(A)$ is the collection of $z\in \C^n$ such that the multiparameter pencil $A(z)=z_1A_1+z_2A_2+\cdots +z_nA_n$ is not invertible. If ${\mathcal B}$ is the group $C^*$-algebra for a discrete group $G$ generated by $A_1,\ A_2,\ ...,\ A_n$ with respect to a representation $\rho$, then $P(A)$ is an invariant of (weak) equivalence for $\rho$. This paper computes the joint spectrum of $R=(1,\ a,\ t)$ for the infinite dihedral group $D_{\infty}=<a,\ t\ |\ a^2=t^2=1>$ with respect to the left regular representation $\lb_{D_{\infty}}$, and gives an in-depth analysis on its properties. A formula for the Fuglede-Kadison determinant of the pencil $R(z)=z_0+z_1a+z_2t$ is obtained, and it is used to compute the first singular homology group of the joint resolvent set $P^c(R)$. The joint spectrum gives new insight into some earlier studies on groups of intermediate growth, through which the corresponding joint spectrum of $(1,\ a,\ t)$ with respect to the Koopman representation $\rho$ (constructed through a self-similar action of $D_{\infty}$ on a binary tree) can be computed. It turns out that the joint spectra with respect to the two representations coincide. Interestingly, this fact leads to a self-similar realization of the group $C^*$-algebra $C^*(D_{\infty})$. This self-similarity of $C^*(D_{\infty})$ is manifested by some dynamical properties of the joint spectrum.
\end{abstract}

\section{Introduction}

The classical spectrum of a linear operator $T$ acting on a Banach space ${\mathcal H}$ is defined as
\[\sigma(T)=\{z\in \C:\ T-zI\ \text{is not invertible}\},\]
and it plays a fundamental role in the study of $T$. Indeed, spectral theory is a centerpiece of operator theory. For a tuple of operators $A=(A_1,\ A_2,\ ...,\ A_n)$, various notions of joint spectrum have been defined to measure joint behavior of the tuple as well as interactions among the elements. If $A$ is a commuting tuple, i.e. $A_iA_j=A_jA_i,\ 1\leq i,j\leq n$, then among others, Taylor spectrum was defined through Koszul complex (cf. \cite{Ho,Ta70}), and it turned out to be a cornerstone in multivariable operator theory. The matter becomes difficult when the tuple is non-commuting, because functional calculus is nearly impossible in this case. Despite the difficulty, some explorations on noncommutative joint spectrum were made, for instance in \cite{Fa}, \cite{Ha73} and \cite{Ta72}.

In many studies of non-commutative tuples, the linear combination (sometimes called multiparameter pencil) \[A(z)=z_1A_1+z_2A_2+\cdots +z_nA_n\]  is a simple and yet pivotal associate of the tuple $A$. It is a generalization of the linear pencil $T-zI$. Invertibility of $A(z)$ is of concern in many areas of mathematics, for instance algebraic geometry, group theory, mathematical physics, PDE, and of course operator theory. We refer the readers to \cite{AJ,At,BG,Sl,Vi,Ya} and the references therein for more information. A case important to this paper is when ${\mathcal B}$ is a group $C^*$ algebra associated with a finitely generated group $G=<g_1,\ g_2,\ \cdots, g_n>$ and a unitary representation $\rho$ on a Hilbert space ${\mathcal H}$. If we let $A_i=\rho(g_i)$, then the invertibility of $A(z)$ reflects the structure of $G$ as well as its representation $\rho$. A good example of this study was made on the following group
\begin{align*}
\G\cong & <a,\ b,\ c,\ d\ |\ 1=a^2=b^2=c^2=d^2=bcd=\sigma^k((ad)^4)\\
&=\sigma^k((adacac)^4),\ k=0,\ 1,\ 2,\ \cdots>, \tag{0.1}
\end{align*}
where in this case
\[\sigma:\ a\to aca,\ b\to d,\ c\to b,\ d\to c\]
is a substitution. It was shown by the first author that $\G$ is of intermediate growth (between polynomial and exponential) (cf. \cite{Gr80}). Then in a series of papers by the first author and his collaborators, the invertibility of the $2$-real variable pencil $Q(\lb,\mu)=-\lb a+b+c+d-(\mu+1)1$ were studied (cf. \cite{BG,Gr84,Gr05}) for the Koopman representation associated with the action on the binary rooted tree and its boundary.

A general study of the invertibility of $A(z)$ was made by the second author in \cite{Ya}, where the notion of projective joint spectrum is introduced as follows.\\

{\bf Definition}. For a tuple $A=(A_1,\ A_2,\ ...,\ A_n)$ of elements in a unital Banach algebra ${\mathcal B}$, its {\em projective joint spectrum} $P(A)$ is the collection of $z\in \C^n$ such that $A(z)=z_1A_1+z_2A_2+\cdots +z_nA_n$ is not invertible in ${\mathcal B}$.\\

A notable distinction of this notion of joint spectrum is that it is ``base free" in the sense that, instead of using $I$ as a base point and looking at the invertibility of $(A_1-z_1I,\ A_2-z_2I,\ ...,\ A_n-z_nI)$ in various constructions, it considers the invertibility of the homogeneous multi-parameter pencil $A(z)$. This feature makes it possible to compute many interesting noncommuting examples.  The projective resolvent set $P^c(A)=\C^n\setminus P(A)$ has some nice properties. For instance, it follows from a general theorem in \cite{ZKKP} that every path-connected component of $P^c(A)$ is a domain of holomorphy. This fact was also proved independently in \cite{Ya} in the case of $C^*$-algebras and in \cite {HY} for general Banach algebras. Further, it was observed in \cite{Ya} that in the case $A$ is commuting $P(A)$ is a union of hyperplanes in $\C^n$, and it is closely related to the Taylor spectrum.

Clearly $0\in \C^n$ is a trivial point in $P(A)$. Because of the homogeneity of $A(z)$,  it is also natural to consider $p(A):=(P(A)\setminus \{0\})/{\C^{\times}}$ which is a bona fide subset in the complex projective space ${\mathbb P}^{n-1}$. As shown in \cite{Ya}, the set $p(A)$ is compact and nonempty for every tuple $A$. Another way to avoid unnecessary complications caused by the trivial point $0$ is to consider
the pencil \[A(z)=I+z_1A_1+z_2A_2+\cdots +z_nA_n,\] in which case $P(A)$ is the collection of $z\in \C^n$ such that $A(z)=I+z_1A_1+z_2A_2+\cdots +z_nA_n$ is not invertible in ${\mathcal B}$. This treatment will be used in several places in the paper to simplify computations, and the meaning of $P(A)$ shall be clear from the context. For more information about the projective joint spectrum, we refer readers to \cite{BCY}, \cite{CSZ}, \cite{CY}, \cite{SYZ} and \cite{Ya}.

Now consider a group $G$ with generators $g_1,\ g_2,\ \cdots,\ g_n$, and let $\rho$ be a unitary representation of $G$ on a Hibert space ${\mathcal H}$. Let $C^*_{\rho}(G)$ denote the $C^*$-algebra generated by $A_i=\rho(g_i)$.
Two representations, say $\rho_1$ on ${\mathcal H}_1$ and $\rho_2$ on ${\mathcal H}_2$ are said to be
equivalent if there is a unitary map $U:\ {\mathcal H}_1\to {\mathcal H}_2$ such that
\[\rho_2(g)=U\rho_1(g)U^{-1},\ \ \forall g\in G.\]
It is obvious that if $\rho_1$ and $\rho_2$ are equivalent representations, then
\[A_{\rho_2}(z)=\sum_iz_i\rho_2(g_i)=U(\sum_iz_i\rho_1(g_i))U^{-1}\] is invertible in $C^*_{\rho_2}(G)$ if and only if $A_{\rho_1}(z)=\sum_iz_i\rho_1(g_i)$ is invertible in $C^*_{\rho_1}(G)$. This indicates that $P(A)$ is an invariant for $\rho$. In fact, $P(A)$ is also invariant under {\em weak equivalence} of unitary representations. This is an important property of projective joint spectrum. We will say more about this in Section 2 and Section 7.

The infinite dihedral group $D_{\infty}=<a,\ t\ |\ a^2=t^2=1>$ plays an important role in group theory. On the one hand, it contains the infinite cyclic normal subgroup $H=<at>$, and $D_{\infty}/H={\mathbb Z}_2$. It is of linear growth and its Cayley graph looks like a line. From this point of view, $D_{\infty}$'s structure is very close to that of ${\mathbb Z}$. But on the other hand, some rather complicated groups can be constructed through $D_{\infty}$. For instance, as demonstrated recently in \cite{Ne}, there is a construction which naturally allows from a minimal action of $D_{\infty}$ to get some complicated nontrivial amenable (and even of intermediate growth) torsion groups. In fact, the main idea of this construction came from a deep analysis of properties of the group $\G$ in (0.1) and of its action on the rooted binary tree and its boundary.

The main theorems in this paper are Theorem 1.1 which gives the projective joint spectrum $P(R)$ for the generating tuple $R=(1,\ a,\ t)$ with respect to the left regular representation $\lb_{D_{\infty}}$, Theorem 4.1 which displays a formula for the Fuglede-Kadison determinant of $R(z)=1+z_1\lb_{D_{\infty}}(a)+z_2\lb_{D_{\infty}}(t)$, Theorem 7.1 which gives a self-similar realization of the group $C^*$-algebra $C^*(D_{\infty})$ through the Koopman representation $\rho$ of $D_{\infty}$ by a measure preserving transformation on the boundary of a binary rooted tree, Theorem 7.2 which gives the joint spectrum with respect to the Koopman representation, and Theorem 8.3 on the dynamics of the joint spectrum $P(R)$ given by a degree-$3$ homogenous polynomial self-map on $\C^3$. This paper is organized as follows.\\

{\bf 0}. Introduction

{\bf 1}.  Joint spectrum through GNS construction

{\bf 2}. Weak containment and maximality

{\bf 3}. Trace of Maurer-Cartan form

{\bf 4}. Fuglede-Kadison determinant

{\bf 5}. On finite dihedral groups

{\bf 6}. Self-similarity

{\bf 7}. Self-similar realization of $C^*(D_{\infty})$

{\bf 8}. Dynamics on joint spectrum

{\bf 9}. Concluding remarks\\

{\bf Acknowledgement.} The authors would like to thank P. de la Harpe, M. Stessin and A. Valette for valuable comments, and P. Kuchment for references. The second author also thanks R. Douglas, Guoliang Yu and the Department of Mathematics at Texas A\&M University for their support and hospitality during his visit.

\section{Joint spectrum through GNS construction}

A state $\phi$ on a $C^*$-algebra ${\mathcal B}$ with unit $I$ is a positive bounded linear functional such that $\phi(I)=1$. And $\phi$ is said to be faithful if $\phi(a^*a)>0$ for every nonzero element $a$. Every faithful state induces an inner product on ${\mathcal B}$ defined by $\langle b,\ a\rangle=\phi(a^*b)$. The completion of ${\mathcal B}$ with respect to the norm induced by the inner product is then a Hilbert space, which we denote by ${\mathcal H}_{\phi}$. Multiplication by elements of ${\mathcal B}$ on ${\mathcal H}_{\phi}$ thus provides a canonical representation of ${\mathcal B}$ on ${\mathcal H}_{\phi}$. If $\phi$ is not faithful, then elements $a$ such that $\phi(a^*a)=0$ form a left ideal $J$ in ${\mathcal B}$, and one can go through the construction of the Hilbert space on ${\mathcal B}/J$. This process is called the Gelfand-Naimark-Segal (GNS) construction. For more details about GNS construction we refer readers to \cite{Av,Da}.

The idea of GNS construction is applicable to group algebras as well. Consider a discrete group  $G$. Its group algebra $\C[G]$ consists of elements of the form \[h=h_01+\sum_{k=1}^{\infty} h_kg_k,\] where $g_k\in G,\ g_k\neq 1$ and $(h_k)$ is a sequence of complex numbers with only a finite number of nonzero values. A conjugate operation, denoted by $^*$, can be defined on $\C[G]$ by
\[h^*=(h_01+\sum_kh_kg_k)^*=\overline{h_0}1+\sum_k\overline{h_k}g_k^{-1},\]
where bar stands for the complex conjugate.
It is a direct computation to check that $(fg)^*=g^*f^*$ for all $f,\ g\in \C[G]$. An element in $\C[G]$ is said to be positive if it is of the form $h^*h$ for some $h\in \C[G]$. Consider the linear functional $tr$ on $\C[G]$ defined by $tr(h)=h_0$. One easily verifies that
\[tr(1)=1,\ \ tr(fg)=tr(gf),\ \ \text{and}\ \ tr(h^*)=\overline{tr(h)},\]
which means $tr$ is a normalized trace on $\C[G]$. Further, since \[tr(h^*h)=|h_0|^2+\sum_k |h_k|^2,\]
$tr(h^*h)=0$ if and only if $h=0$, meaning that $tr$ is faithful. The GNS construction then gives rise to the Hilbert space ${\mathcal H}_{tr}$, as well as the canonical respresentation $\lb_{G}$ such that
\[\lb_G(g)h=gh=h_0g+\sum_{k=1}^{\infty}h_kgg_k,\ \ g\in G,\ h\in {\mathcal H}_{tr}.\]
Further, for any $h,\ h'\in {\mathcal H}_{tr}$
\[\langle gh',\ gh\rangle=tr((gh)^*(gh'))=tr(h^*g^*gh')=tr(h^*g^{-1}gh')=tr(h^*h')=\langle h',\ h\rangle,\]
so $\lb_G$ is a unitary representation.

${\mathcal H}_{tr}$ can be canonically identified with $l^2(G)$ by the unitary map $U(g)=\delta_g,\ g\in G$, where $\delta_g$ takes value $1$ at $g$ and $0$ at other elements in $G$.
One checks that $\lb_G$ is equivalent to the left regular representation of $G$ on $l^2(G)$. Therefore we deal with the $C^*$-algebra $C^*_{\lb_G}(G)$ generated by left regular representation, which in the case when $G$ is amenable coincides with the full $C^*$-algebra of $G$ denoted by $C^*(G)$.

The representation $\lb_G$ gives rise to a natural extension of the trace $tr$ from $\C[G]$ to the von Neumann algebra $L(G)$ generated by $\lb_G(g),\ g\in G$. Noting that $1\in {\mathcal H}_{tr}$, for $F\in L(G)$ we can define
\begin{equation*}
trF=\langle F 1,\ 1\rangle.\tag{1.0}
\end{equation*}

The infinite dihedral group $D_{\infty}$ is isomorphic to the free product ${\mathbb Z}_2\ast {\mathbb Z}_2$, and as remarked earlier, it contains the normal subgroup $H=<at>$ which is an infinite cyclic group with $D_{\infty}/H={\mathbb Z}_2$ (this, in particular, implies that $D_{\infty}$ is amenable). So elements in $D_{\infty}$ are of the form $(at)^k$ or $t(at)^k$, where $k\in {\mathbb Z}$, for example, $a=t(at)^{-1}$. Hence the
complex group algebra $\C[D_{\infty}]$ consists of elements of the form
\[h=\sum_kb_k(at)^k+\sum_j c_jt(at)^j,\]
where $\{b_k, c_j: k,j\in {\mathbb Z}\}\subset \C$ has only a finite number of nonzero elements.

As defined above, $tr$ is a faithful state on $\C[D_{\infty}]$. When restricted to the subalgebra $\C[H]$, the GNS construction gives rise to the Hilbert space
\[L=\{f= \sum_{-\infty}^{\infty} f_j(at)^j:\ \sum_{-\infty}^{\infty} |f_j|^2<\infty\},\]
with inner product
\[\langle g,\ f\rangle=tr(f^*g)=\sum_{-\infty}^{\infty}\overline{f_j}g_j.\]
Observe that$ \{(at)^j:\ j\in {\mathbb Z}\}$ is an orthonormal basis for $L$. Hence multiplication by $at$, which we denote by $T$, is unitarily equivalent to the classical bilateral shift operator: multiplication by $e^{i\theta}$ on $L^2({\mathbb T},\ \frac{d\theta}{2\pi})$.

Since $\C[D_{\infty}]=\C[H]\oplus t\C[H]$, the GNS construction for $\C[D_{\infty}]$ then gives the Hilbert space orthogonal direct sum ${\mathcal H}_{tr}=L\oplus tL$. On $L\oplus tL$, one checks that
\[a(at)^j=t(at)^{j-1}=tT^*((at)^j)\in tL,\ \text{and}\ \ a[t(at^j)]=(at)^{j+1}=T((at)^j)\in L.\] Moreover,
$t(at)^j\in tL$ and $t[t(at)^j]=(at)^j$. Since $\lb_{D_{\infty}}(a)$ and $\lb_{D_{\infty}}(t)$ are multiplications by $a$ and respectively $t$ on $L\oplus tL$, we have
\begin{equation*}
\lb_{D_{\infty}}(a)=\left(
      \begin{array}{cc}
        0 & Tt  \\
        tT^* & 0 \\
      \end{array}
    \right),\quad
\lb_{D_{\infty}}(t)=\left(
      \begin{array}{cc}
        0 & t  \\
        t & 0 \\
      \end{array}
    \right).
\end{equation*}
Now consider the map $S: L\oplus tL\to L\oplus L$ defined by
\begin{equation*}
S=\left(
      \begin{array}{cc}
        I_L & 0 \\
        0 & t \\
      \end{array}
    \right),
\end{equation*}
where $I_L$ is the identity operator on $L$. Clearly, $S$ is unitary and $S^2=I$. Further, one verifies that
\begin{equation*}
\lb_{D_{\infty}}(a)=S\left(
      \begin{array}{cc}
        0 & T  \\
        T^* & 0 \\
      \end{array}
    \right)S,\quad
\lb_{D_{\infty}}(t)=S\left(
      \begin{array}{cc}
        0 & I_L  \\
        I_L & 0 \\
      \end{array}
    \right)S.\tag{1.1'}
\end{equation*}
So, up to the unitary $S$, we can write
\begin{equation}
\lb_{D_{\infty}}(a)=\left(
      \begin{array}{cc}
        0 & T  \\
        T^* & 0 \\
      \end{array}
    \right),\quad
\lb_{D_{\infty}}(t)= \left(
      \begin{array}{cc}
        0 & I_L  \\
        I_L & 0 \\
      \end{array}
    \right).
\end{equation}
The representation $\lb_{D_{\infty}}$ in (1.1) extends to a representation of $D_{\infty}$ on $L\oplus L$. The $C^*$ algebra generated by $\lb_{D_{\infty}}(a)$ and $\lb_{D_{\infty}}(t)$ through this representation is then isometrically isomorphic to $C^*(D_{\infty})$.

The representation $\lb_{D_{\infty}}$ in (1.1) enables us to compute the projective joint spectrum in concern. For the tuple $R=(1,\ a,\ t)$, we let $R_{\lb_{D_{\infty}}}(z)=z_0\lb_{D_{\infty}}(1)+z_1\lb_{D_{\infty}}(a)+z_2\lb_{D_{\infty}}(t)$. Sections 1-5 are mostly concerned with the representation $\lb_{D_{\infty}}$ in (1.1), so for simplicity we write $R_{\lb_{D_{\infty}}}(z)$ as $R(z)=z_0+z_1a+z_2t$. Then $P(R)$ is the set of $z\in \C^3$ such that $R(z)$ is not invertible in $C^*(D_{\infty})$.

\begin{thm}
For $R(z)=z_0+z_1a+z_2t$, with respect to the representation $\lb_{D_{\infty}}$ the projective joint spectrum
\[P(R)=\bigcup_{-1\leq x\leq 1}\{z\in \C^3:\ z_0^2-z_1^2-z_2^2-2z_1z_2x=0\}.\]
\end{thm}

\pf  It is well-known that the bilateral shift $T$ is a unitary with spectrum $\sigma(T)$ equal to the unit circle ${\mathbb T}$, and its spectral resolution is given by
\begin{equation}
T=\int_{\mathbb T}\lb dE(\lb).
\end{equation}
where $E(\lb)$ is the associated projection-valued spectral measure. Using (1.1), we can now write
\begin{eqnarray}
R(z)=z_0+z_1a+z_2t=\left(
      \begin{array}{cc}
        z_0 & z_1T+z_2  \\
        z_1T^*+z_2 & z_0 \\
      \end{array}
    \right).
\end{eqnarray}

We consider two cases.

If $z_0=0$, then $R(z)$ is invertible if and only if both $z_1T+z_2$ and $z_1T^*+z_2$ are invertible, or by (1.2)
\[(z_1\lb +z_2)(z_1\bar{\lb}+z_2)=z_1^2+z_2^2+z_1z_2(\lb+\bar{\lb})\neq 0\] for every $\lb\in {\mathbb T}$.

If $z_0\neq 0$, then by factorization we may consider the case $z_0=1$, i.e.  $R(z)=1+z_1a+z_2t$. Recall that for a block matrix,
\begin{equation}
\begin{bmatrix}
A & B \\ C & D
 \end{bmatrix}^{-1}
=
\begin{bmatrix}
A^{-1}+A^{-1}BK^{-1}CA^{-1} & -A^{-1}BK^{-1} \\
-K^{-1}CA^{-1} & K^{-1}
\end{bmatrix},
\end{equation}
when both $A$ and the Schur complement $K=D-CA^{-1}B$ are invertible. Further, when $A$ is invertible, the block matrix is invertible if and only if $D-CA^{-1}B$ is invertible (\cite{GN07,LS}).
So in this case $R(z)$ is invertible if and only if $1-(z_1T+z_2)(z_1T^*+z_2)$ is invertible on $L$, or by (1.2)
\[1-(z_1\lb+z_2)(z_1\bar{\lb}+z_2)=1-z_1^2-z_2^2-z_1z_2(\lb+\bar{\lb})\neq 0\]
for every $\lb\in \td$.

Summarizing these two cases, we have $R(z)$ is not invertible if and only if
\[\ z_0^2-z_1^2-z_2^2-z_1z_2(\lb+\bar{\lb})=z_0^2-z_1^2-z_2^2-2z_1z_2\cos \theta =0,\]
for some $\lb=e^{i\theta}\in \td$. Setting $\cos \theta=x$, we have the theorem.
\zb

We end this section by an observation that will be used later.
\begin{corr}
For $R(z)=z_0+z_1a+z_2t$, with respect to the representation $\lb_{D_{\infty}}$ the joint resolvent set $P^c(R)$ is path connected.
\end{corr}
\pf The joint spectrum $P(R)$ is displayed in Theorem 1.1. We first look at the case $z_0\neq 0$. Since $R(z)$ is homogenous in $z$ and we can write $(z_0,\ z_1,\ z_2)=z_0(1,\ z_1/z_0,\ z_2/z_0)$ in this case, we shall prove without loss of generality with the assumption $z_0=1$. To avoid possible confusion, in here we let $R_{*}(z)=1+z_1a+z_2t$. Then by Theorem 1.1 we have
\[P(R_{*})=\bigcup _{-1\leq x\leq 1}\{z\in \C^2:\ 1-z_1^2-z_2^2-2z_1z_2x=0\}.\]
We now check that $P^c(R_{*})$ is path-connected. To this end, we check that every point $\lb=(\lb_1,\ \lb_2)\in P^c(R_{*})$ is path connected to $(0,\ 0)$. By possibly choosing a point in a small ball centered at $\lb$, and using symmetry of $P(R_{*})$, we may assume without loss of generality that $|\lb_1|>|\lb_2|>0$. Consider the complex line \[C_{\lb}=\C\lb=\{(w\lb_1,\ w\lb_2): w\in \C\}.\] Clearly,
$C_{\lb}\cap P^c(R_{*})$ is a subset of $P^c(R_{*})$ and contains $(0,\ 0)$ and $(\lb_1,\ \lb_2)$. So it is sufficient to check that $C_{\lb}\cap P^c(R_{*})$ is path connected. Observe that
\[C_{\lb}\cap P(R_{*})=\bigcup_{-1\leq x\leq 1}\{w\in \C:\ 1-w^2(\lb_1^2+\lb_2^2)-2w^2\lb_1\lb_2x=0\}.\]
Solving for $w^2$, we have $w^2=(\lb_1^2+\lb_2^2+2\lb_1\lb_2x)^{-1}$ and hence we have two solution curves \[w_{\pm}(x)=\pm (\lb_1^2+\lb_2^2+2\lb_1\lb_2x)^{-1/2},\ \ -1\leq x\leq 1.\]
Clearly $w_{\pm}(-1)=\pm 1/(\lb_1-\lb_2)$ and $w_{\pm}(1)=\pm 1/(\lb_1+\lb_2)$. All four points are distinct.  So as $x$ moves from $-1$ to $1$ along the real line, the values of $w_+(x)$ form a simple path connecting $1/(\lb_1-\lb_2)$ to $1/(\lb_1+\lb_2)$; and the values of $w_{-}(x)$ form a simple path connecting $-1/(\lb_1-\lb_2)$ to $-1/(\lb_1+\lb_2)$. Now we verify the following two facts.

1. The two paths are not intersecting. This is because if there are $x_1, x_2\in [-1,\ 1]$ such that $w_+(x_1)=w_{-}(x_2)$, then $w^2_+(x_1)=w^2_{-}(x_2)$, which implies
\[\lb_1^2+\lb_2^2+2\lb_1\lb_2x_1=\lb_1^2+\lb_2^2+2\lb_1\lb_2x_2\]
 and hence $x_1=x_2$.  It follows that $w_+(x_1)=w_{-}(x_1)$ which happens only when both are $0$, which is impossible.

2. Neither curve self-intersects. If there are $x_1, x_2\in [-1,\ 1]$ such that $w_+(x_1)=w_{+}(x_2)$ (or $w_-(x_1)=w_{-}(x_2)$), then by similar arguments, we will have $x_1=x_2$.

In summary, the set $C_{\lb}\cap P(R_{*})$ is a disjoint union of two simple closed curves. Hence its complement in $C_{\lb}$ i.e. $C_{\lb}\cap P^c(R_{*})$ is path connected. This concludes that $P^c(R_{*})$, and hence $P^c(R)\cap \{z_0\neq 0\}$, is path connected.

For the case $z_0=0$, we see by Theorem 1.1 that a fixed $(0,\ z_1,\ z_2)\in P^c(R)$ if and only if $z_1^2+z_2^2+2z_1z_2x\neq 0$ for all $x\in [-1,\ 1]$. Since $[-1,\ 1]$ is compact, there is a constant $\delta$ such that
\[ |z_1^2+z_2^2+2z_1z_2x|\geq \delta>0,\ \  \forall \ x\in [-1,\ 1].\]
Pick $z_0$ such that $0<|z_0|^2<\delta$, then the path $(tz_0,\ z_1,\ z_2),\ 0\leq t\leq 1$ lies in $P^c(A)$ by Theorem 1.1, and it connects $(0,\ z_1,\ z_2)$ to $(z_0,\ z_1,\ z_2)$, which in turn is connected to $(1,\ 0,\ 0)$.

\zb

\section{Weak containment and Maximality}

As we have remarked in Section 0, projective joint spectrum is an invariant for representations of groups
up to unitary equivalence. We begin this section with a simple example.
\begin{exam}
For a fixed $\theta\in (0,\ \pi)$, we consider the following two-dimensional irreducible representation $\rho_{\theta}$ of $D_{\infty}$ given by
\[\rho_{\theta}(a)=\begin{bmatrix}
         0 & e^{i \theta}\\
        e^{-i \theta} & 0
        \end{bmatrix},\ \ \ \rho_{\theta}(t)=\begin{bmatrix}
                                          0 & 1\\
                                          1 & 0
                                          \end{bmatrix}.\]

It is known (\cite{Hal,RS}) that every unitary irreducible representation of $D_{\infty}$ is either one dimensional or of the form $\rho_{\theta}$.
Clearly, the classical spectra $\sigma(\rho_{\theta}(a))=\sigma(\rho_{\theta}(t))=\{\pm 1\}$, and their characters $Tr\rho_{\theta}(a)=Tr\rho_{\theta}(t)=0$, none of which reflects the representation $\rho_{\theta}$. But if we consider the pencil $R_{\theta}(z)=z_0+z_1\rho_{\theta}(a)+z_2\rho_{\theta}(t)$, then \[det R_{\theta}(z)=z_0^2-z_1^2-z_2^2-2z_1z_2\cos (\theta),\]
and hence the joint spectrum $P(R_{\theta})$ is the variety $\{det R_{\theta}(z)=0\}$. Therefore, if $\cos (\theta_1)\neq \cos (\theta_2)$ then the two representations $\rho_{\theta_1}$ and $\rho_{\theta_2}$ are not equivalent.
\end{exam}
Observe that the variety $\{det R_{\theta}(z)=0\}$ is a slice in $P(R)$ in Theorem 1.1. Since every unitary representation is a direct integral of irreducible ones, Theorem 1.1 makes one wonder if $P(R)$ in Theorem 1.1 is maximal in some sense. This is indeed so, and it is a consequence of the maximality of left regular representation for amenable groups in the relation of weak containment.\\

{\bf Definition.}  Consider two representations $\pi$ and $\rho$ of a discrete group $G$ in Hilbert spaces ${\mathcal H}$ and ${\mathcal K}$, respectively. One says that $\pi$ is weakly contained in $\rho$ (denoted by $\pi\prec \rho$) if for every $x\in {\mathcal H}$, every finite subset $F\subset G$, and every $\epsilon>0$, there exist $y_1,\ y_2,\ \cdots,\ y_n$ in ${\mathcal K}$ such that for all $g\in F$
\[|\langle \pi(g)x,\ x\rangle-\sum_{i=1}^{n}\langle \rho(g)y_i,\ y_i\rangle|<\epsilon.\]

It is known that if $\pi\prec \rho$ then the map $\rho(m)\to \pi(m),\ m\in \C[G]$ extends to a unital homomorphism from $C^*_{\rho}(G)$ onto $C^*_{\pi}(G)$. In particular, this implies that if $\rho(m)$ is invertible in $C^*_{\rho}(G)$ then $\pi(m)$ is invertible in $C^*_{\pi}(G)$. So for $\{g_1,\ g_2,\ \cdots,\ g_n\}\subset G$ if we let \[A_{\rho}=(\rho(g_1),\ \rho(g_2), \cdots,\ \rho(g_n)),\] and $A_{\pi}$ likewise, then $\pi\prec \rho$ implies $P(A_{\pi})\subset P(A_{\rho})$.  If $\pi\prec \rho$ and $\rho \prec \pi$, then $\pi$ and $\rho$ are said to be weakly equivalent, and we denote this by $\pi\sim\rho$. Clearly, we have $P(A_{\pi})= P(A_{\rho})$ in this case. These facts indicate that projective joint spectrum provides a measurement for weak containment.
The  following  theorem  of  Hulanicki and  Reiter (cf. \cite{BHV} Theorem G 3.2) shows that for amenable groups regular representation is maximal  with respect to the relation given by the weak containment.\\

\noindent {\bf Theorem}. Let $G$ be a locally compact group. The following properties are equivalent:

(i) $G$ is amenable;

(ii) ${\bf 1}_G \prec \lb_G$;

(iii) $\pi \prec \lb_G$ for every unitary representation $\pi$ of $G$.\\

Here ${\bf 1}_G$ stands for the trivial representation. So for locally compact amenable groups, $\lb_G$ is maximal among unitary representations with respect to weak containment. These observations lead to the following

\begin{prop}
Let $G$ be a locally compact amenable group and $\{g_1,\ g_2,\ \cdots,\ g_n\}$ be any finite subset in $G$. Then

(i) for every unitary representation $\pi$ of $G$, we have $P(A_{\pi})\subset P(A_{\lb _G})$.

(ii) $P(A_{\lb _G})$ contains the hyperplane $\{z_1+z_2+\cdots +z_n=0\}$.
\end{prop}
The second statement follows from the fact ${\bf 1}_G\prec \lb_G$ and \[R_{{\bf 1}_G}(z)=(z_1+z_2+\cdots +z_n)I.\]
A special case of (ii) in Proposition 2.2 is when $n=2$ with $g_1=1$ and $g_2=g$ an arbitrary element in $G$. Statement (ii) then implies that $I-\lb_G(g)$ is not invertible, or equivalently $1\in \sigma(\lb_G(g))$ for every $g\in G$.

So regarding $D_{\infty}$, an immediate consequence of Proposition 2.2 is that if $\pi$ is a unitary representation of $D_{\infty}$, and $R_{\pi}=z_01+z_1\pi(a)+z_2\pi(t)$, then
$P(R_{\pi})$ is a subset of the joint spectrum $P(R)$ in Theorem 1.1.

Now we turn to $p=(1-a)/2$ and $q=(1-t)/2$. One sees that $p$ and $q$ are two projections in $C^*(D_{\infty})$. They are called projections in ``generic position" in \cite{RS}. To be precise, if $p$ and $q$ are two projections in a unital $C^*$-algebra ${\mathcal B}$ such that the $C^*$-subalgebra generated by $I,\ p$ and $q$ is isomorphic to $C^*(D_{\infty})$ then $p$ and $q$ are said to be in generic position.

Let $A=(I,\ p,\ q)$ and $A(z)=z_0+z_1p+z_2q$. One checks easily that
\[A(z)=\big(z_0+\frac{z_1+z_2}{2}\big)-\frac{z_1}{2}a-\frac{z_2}{2}t.\]
So by Theorem 1.1 $A(z)$ is not invertible if and only if
\begin{align*}
&\big(z_0+\frac{z_1+z_2}{2}\big)^2-\big(\frac{z_1}{2}\big)^2-\big(\frac{z_2}{2}\big)^2-\frac{z_1z_2}{2}\cos \theta\\
&=z_0^2+z_0(z_1+z_2)+\frac{z_1z_2}{2}(1-\cos\theta)\\
&=z_0^2+z_0(z_1+z_2)+z_1z_2\sin^2\frac{\theta}{2}=0,
\end{align*}
for some $0\leq \theta \leq 2\pi$. Setting $x=\sin^2\frac{\theta}{2}$ we have
\[P(A)=\bigcup_{0\leq x\leq 1}\{z\in \C^3:\ z_0^2+z_0(z_1+z_2)+z_1z_2x=0\}.\]
The $p$ and $q$ are universal in the following sense (cf. \cite{RS}): if $p',\ q'$ is an arbitary pair of projections in a unital $C^*$-algebra $B$, then there is a unital homomorphism $\phi:\ C^*(D_{\infty})\longrightarrow B$ such that $\phi(p)=p',\ \phi(q)=q'$. Therefore, if $A(z)$ is invertible in $C^*(D_{\infty})$, then
\[\phi(A(z))=z_0\phi(I)+z_1\phi(p)+z_2\phi(q)=z_0I+z_1p'+z_2q':=A'(z)\]
is invertible in $B$. Summarizing these facts we have
\begin{corr}
If $p$ and $q$ are projections in generic position and $A=(I,\ p,\ q)$, then
\[P(A)=\bigcup_{0\leq x\leq 1}\{z\in \C^3:\ z_0^2+z_0(z_1+z_2)+z_1z_2x=0\}.\]
And for an arbitrary pair of projections $p',\ q'$, $P(A')$ is a subset of $P(A)$.
\end{corr}

Two more observations are worth mentioning.

{\bf 1.} Considering the case $x=0$ in Corollary 2.3, we see that $P(A)$ contains the slice \[\{z\in \C^3:\ z_0(z_0+z_1+z_2)=0\},\] i.e. the hyperplanes $\{z_0+z_1+z_2=0\}$ and $\{z_0=0\}$ are in $P(A)$. In the case $z_0=0$, $A(z)=z_1p+z_2q$ with $z_1$ and $z_2$ arbitrary. So it implies that for projections $p$ and $q$ the linear combination $z_1p+z_2q$ is not invertible for any $z_1$ and $z_2$. We state this fact as
\begin{corr}
If $p$ and $q$ are projections in generic position, then
$z_1p+z_2q$ is not invertible for any complex numbers $z_1$ and $z_2$.
\end{corr}

Now consider the group
\[\tilde{D}=<a_1,a_2,a_3|\ a_1^2=a_2^2=a_3^2=1>\cong {\mathbb Z}_2\ast {\mathbb Z}_2 \ast {\mathbb Z}_2.\]
Since it contains the free subgroup on two generators $<a_1^{-1}a_2, a_1^{-1}a_3>\cong {\mathbb Z}\ast {\mathbb Z}$, $\tilde{D}
$ is not amenable and the full group $C^*$-algebra $C^*(\tilde{D})$ has a rather complicated structure (cf. \cite{RS}). In particular, Proposition 2.2 is not valid in this case. It is thus a natural question whether or not one can compute the joint spectrum for the tuple $(1, a_1, a_2, a_3)$, or equivalently for the tuple $(I, q_1, q_2, q_3)$ where $q_i=(1-a_i)/2,\ i=1,2,3$.\\

\noindent {\bf Problem.} Determine $z\in \C^4$ such that $z_0 I+z_1q_1+z_2q_2+z_3q_3$ is not invertible.\\

{\bf 2.} In \cite{Pe} Pedersen displayed a seemingly different representation of $C^*(D_{\infty})$ in terms of matrix-valued functions: there is an isomorphism of $C^*(D_{\infty})$ onto
\[{\mathcal A}=\{f\in C([0,\ 1], M_2(\C)):\ f(0),\ f(1)\ \text{are diagonal}\},\]
which carries $p=(1-a)/2$ and $q=(1-t)/2$ into the functions
\[p(x)=\begin{pmatrix} x & \sqrt{x(1-x)} \\ \sqrt{x(1-x)} & 1-x \end{pmatrix}, \
  q(x)=\begin{pmatrix} 1 & 0 \\ 0 & 0 \end{pmatrix}.\]
Here $C([0,\ 1], M_2(\C))$ is the set of $2\times 2$ complex matrix-valued continuous functions on $[0,\ 1]$.
Under this isomorphism, $a$ and $t$ are represented by
\[a(x)=\begin{pmatrix} 1-2x & -2\sqrt{x(1-x)} \\ -2\sqrt{x(1-x)} & 2x-1 \end{pmatrix},\
t(x)=\begin{pmatrix} -1 & 0 \\ 0 & 1 \end{pmatrix}.\]
The joint spectrum of $(1,\ a(x),\ t(x))$ can be computed directly, and it turns out to coincide with $P(R)$ in Theorem 1.1. As a matter of fact, if we let $V$ be the unitary matrix
\[\frac{1}{\sqrt{2}}\begin{pmatrix} 1 & i \\ -1 & i \end{pmatrix},\]
then from (1.1) we have that
\[V^*tV=\begin{pmatrix} -1 & 0 \\ 0 & 1 \end{pmatrix}.\]
Using the fact that $T$ is unitarily equivalent to multiplication by $e^{i\theta}$ on $L^2({\mathbb T})$, we also have
\[V^*aV=\frac{1}{2}\begin{pmatrix} -(T+T^*) & i(T-T^*) \\ i(T-T^*) & T+T^* \end{pmatrix}
\cong \begin{pmatrix} -\cos \theta &  -\sin\theta \\
                           -\sin \theta &  \cos{\theta} \end{pmatrix}.\]
Setting $x= \cos^2\frac{\theta}{2}$ and using trigonometric identities, we have Pedersen's representation of $C^*(D_{\infty})$. So in fact Pedersen's representation and the left regular representation $\lb_{D_{\infty}}$ are unitarily equivalent.

\section{Trace of Maurer-Cartan form}

For a tuple $A=(A_1,\ A_2,\ \cdots,\ A_n)$ of elements in a unital Banach algebra ${\mathcal B}$, recall that $A(z)=z_1A_1+z_2A_2+\cdots +z_nA_n$. The Maurer-Cartan type ${\mathcal B}$-valued $1$-form on $P^c(A)$ defined by
\[\omega_A(z)=A^{-1}(z)dA(z)=\sum_{j=1}^nA^{-1}(z)A_jdz_j,\ \ z\in P^c(A),\]
is an important subject of study in \cite{Ya}, and it is shown to contain much information about the topology of $P^c(A)$ which can be read by invariant linear functionals or by cyclic cocycles (cf. \cite{CY,Ya}). For example, it is indicated in \cite{Ya} that if ${\mathcal B}$ possesses a trace $tr$, then
$tr\omega_A(z)$ is a nontrivial element in the de Rham cohomology $H^1(P^c(A),\C)$.

Representation $\lb_{D_{\infty}}$ enables us to do a more in-depth study of $R(z)=z_0+z_1a+z_2t$. In this section we shall compute $trR^{-1}(z)$ and $tr\omega_R$. This computation is important for our discussion on the Fuglede-Kadison determinant in the next section.  For simplicity, we shall do the computation for the case $z_0=1$, and the general case is a natural extension.

We first compute $trR^{-1}(z)$. Letting $K(z)=1-(z_1T+z_2)(z_1T^*+z_2)$ and using (1.4), we have that on $L\oplus L$,
\begin{equation*}
R^{-1}(z)=\begin{bmatrix}
 1+(z_1T+z_2)K^{-1}(z)(z_1T^*+z_2) & -(z_1T+z_2)K^{-1}(z) \\
-K^{-1}(z)(z_1T^*+z_2) & K^{-1} (z)
\end{bmatrix}.\tag{3.0}
\end{equation*}
Note that through the unitary $S$ in (1.1'), the entries in $R^{-1}(z)$ above can be identified with elements in the von Neuman algebra $L(D_{\infty})$. The trace $tr$ on $L(D_{\infty})$ defined as in (1.0) can be naturally extended to $2\times 2$ matrices with $L(D_{\infty})$ entries by the definition
\[tr \begin{pmatrix} a_{11} & a_{12} \\ a_{21} & a_{22} \end{pmatrix}:=\frac{1}{2}tr(a_{11}+a_{22}).\]
Note that $tr(I)=tr(I_L)=1$. Then using (1.2), we have
\begin{align}
trR^{-1}(z)&=\frac{1}{2}tr\big( 1+(z_1T+z_2)K^{-1}(z)(z_1T^*+z_2)+K^{-1} (z)\big)\\
&=\frac{1}{2}tr\int_{\td}1+\frac{(z_1\lb+z_2)(z_1\bar{\lb}+z_2)+1}{1-z_1^2-z_2^2-z_1z_2(\lb+\bar{\lb})}dE(\lb)\\
&=\frac{1}{2}\int_{\td}\frac{2trdE(\lb)}{1-z_1^2-z_2^2-z_1z_2(\lb+\bar{\lb})}.
\end{align}
Now we take a closer look at $tr$. In Section 1, the GNS construction on $\C[H]$ produced the Hilbert space $L$ with inner product $\langle g,\ f\rangle=tr(f^*g)$. Now consider the linear map $U:\ L\to L^2({\mathbb T},\ \frac{d\theta}{2\pi})$ defined by $U((at)^n)=e^{ni\theta}$, $n\in {\mathbb Z}$. One easily checks that $U$ is a unitary. In particular,
\begin{align}
tr(f^*g)=\langle g,\ f\rangle=\int_0^{2\pi}\overline{(Uf)(e^{i\theta})}(Ug)(e^{i\theta})\frac{d\theta}{2\pi}.
\end{align}
Further, for any $f$ and $g$ in $\C[H]$, functional calculus gives $f=(Uf)(T)$ and $g=(Ug)(T)$. Hence by (1.2),
\[tr(f^*g)=tr\int_{\mathbb T}\overline{(Uf)(\lb)}(Ug)(\lb)dE(\lb)=\int_0^{2\pi}\overline{(Uf)(e^{i\theta})}(Ug)(e^{i\theta})tr(dE(e^{i\theta})).\]
Comparing this with (3.4), we have
\begin{equation}
tr (dE(e^{i\theta}))=dsrE(e^{i\theta})=\frac{d\theta}{2\pi}.
\end{equation}
So it follows from (3.3) that
\begin{equation}
trR^{-1}(z)=\frac{1}{2\pi}\int_{0}^{2\pi}\frac{d\theta}{1-z_1^2-z_2^2-2z_1z_2\cos\theta}.
\end{equation}
We thus obtain the following
\begin{prop}
For $R(z)=I+z_1a+z_2t$, we have
\[trR^{-1}(z)=\frac{1}{2\pi}\int_{0}^{2\pi}\frac{d\theta}{1-z_1^2-z_2^2-2z_1z_2\cos\theta},\ \ z\in P^c(R).\]
\end{prop}
Observe that $trR^{-1}(z)$ is holomorphic on $P^c(R)$ and can not be extended holomorphically into a neighborhood of any points in $P(R)$, indicating that $P^c(R)$ is a domain of holomorphy as remarked in Introduction.

Now for $R(z)=1+z_1a+z_2t$, we consider the Maurer-Cartan form
\[\omega_R(z)=R^{-1}(z) dR(z)=R^{-1}(z)(adz_1+tdz_2),\]
and compute $tr\omega_R(z)$. Writing the matrix form of $R^{-1}(z)$ (cf. (3.10)) as $(R^{jk}(z))_{2\times 2}$, we check that
\begin{align*}
\omega_R(z)&=R^{-1}(z)\left(\begin{bmatrix} 0 & T \\ T^* & 0 \end{bmatrix} dz_1
                                             +\begin{bmatrix} 0 & 1 \\ 1 & 0 \end{bmatrix} dz_2    \right)\\
&=  \begin{bmatrix}
        R^{12}T^* & R^{11}T \\
        R^{22}T^* & R^{21}T
      \end{bmatrix} dz_1+  \begin{bmatrix}
        R^{12} & R^{11} \\
        R^{22} & R^{21}
      \end{bmatrix}dz_2.
\end{align*}
Hence by (1.2), (3.5), and using computations similar to (3.1) we have
\begin{align*}
tr \omega_R(z)&=\frac{1}{2}tr(R^{12}T^*+R^{21}T)dz_1+\frac{1}{2}tr(R^{12}+R^{21})dz_2\\
&=\frac{-1}{2}\int_{\td}\frac{(z_1\lb+z_2)\bar{\lb}+(z_1\bar{\lb}+z_2)\lb}{1-z_1^2-z_2^2-        z_1z_2(\lb+\bar{\lb})}dsrE(\lb)dz_1\\
&+\frac{-1}{2}\int_{\td}\frac{(z_1\lb+z_2)+(z_1\bar{\lb}+z_2)}{1-z_1^2-z_2^2-z_1z_2(\lb+\bar{\lb})}dsrE(\lb)dz_2\\
&=\frac{-1}{2}\int_{\td}\frac{(2z_1+z_2(\lb+\bar{\lb}))dz_1+(z_1(\bar{\lb}+\lb)+2z_2)dz_2}{1-z_1^2-z_2^2-z_1z_2(\lb+\bar{\lb})}dsrE(\lb)\\
&=\frac{1}{4\pi}\int_0^{2\pi}\partial \log (1-z_1^2-z_2^2-2z_1z_2\cos\theta)d\theta\\
&=\partial \big(\frac{1}{4\pi}\int_0^{2\pi} \log (1-z_1^2-z_2^2-2z_1z_2\cos\theta)d\theta\big),
\end{align*}
where $\partial f=\frac{\partial f}{\partial z_1}dz_1+\frac{\partial f}{\partial z_2}dz_2$.
We summarize this calculation as
\begin{prop}
For $R(z)=I+z_1a+z_2t$,
\[tr\omega_R(z)=\partial \big(\frac{1}{4\pi}\int_0^{2\pi} \log (1-z_1^2-z_2^2-2z_1z_2\cos\theta)d\theta\big).\]
\end{prop}
Three remarks are in order.

{\bf 1}. Clearly $tr\omega_R(z)$ is a closed $1$-form. Since  $\log (1-z_1^2-z_2^2-2z_1z_2\cos\theta)$ is not globally defined on $P^c(R)$, $tr\omega_R(z)$ is not exact and hence is a nontrivial element in $H^1(P^c(R),\C)$. This fact will become more apparent when we couple it with the fundamental group of $P^c(R)$ in the next section (cf. (4.2) and Corollary 4.4). Further, if $\gamma=\{z(s): 0\leq s\leq 1\}$ is a piece-wise smooth path, then when restricted to $\gamma$ we have
\begin{equation}
tr\omega_R(z(s))= \frac{d}{ds}\big(\frac{1}{4\pi}\int_0^{2\pi} \log (1-z_1^2(s)-z_2^2(s)-2z_1(s)z_2(s)\cos\theta)d\theta\big)ds.
\end{equation}
This equality will be used in (4.2) as well.

{\bf 2}. Observe further that for any fixed $z\in P(R)$ such that $z_1z_2\neq 0$, $\log (1-z_1^2-z_2^2-2z_1z_2\cos\theta)$ is integrable with respect to $\theta$ over $[0,\ 2\pi]$. Hence the function
\[ \frac{1}{4\pi}\int_0^{2\pi} \log (1-z_1^2-z_2^2-2z_1z_2\cos\theta)d\theta\]
can be extended holomorphically into this part of the spectrum, even though it is not globally defined. We shall have more to say about this as well in the next section.

{\bf 3}. If we consider $R(z)=z_0I+z_1a+z_2t$, then with just a bit more computation one can check that
\[trR^{-1}(z)=\frac{1}{2\pi}\int_{0}^{2\pi}\frac{z_0d\theta}{z_0^2-z_1^2-z_2^2-2z_1z_2\cos\theta},\ \ z\in P^c(R),\]
and
\[tr\omega_R(z)=\partial \big(\frac{1}{4\pi}\int_0^{2\pi} \log (z_0^2-z_1^2-z_2^2-2z_1z_2\cos\theta)d\theta\big),\]
where $\partial f=\frac{\partial f}{\partial z_0}dz_0+\frac{\partial f}{\partial z_1}dz_1+\frac{\partial f}{\partial z_2}dz_2$.

\section{Fuglede-Kadison determinant}

For a unital Banach algebra $ {\mathcal B}$, we let $GL({\mathcal B})$ denote the open subset of invertible elements in ${\mathcal B}$. In \cite{FK}, Fuglede and Kadison defined the notion of determinant for invertible elements $x$ in a finite von Neumann algebra $ {\mathcal B}$ with a normalized trace $tr$ as \[\det x=\exp(tr\log \sqrt{x^*x}),\]
and made a general study on its properties. In particular, they showed that $\det$ is a homomorphism from $GL({\mathcal B})$ to the multiplicative group ${\mathbb R}_+$ of positive real numbers. The FK-determinant can be extended analytically to non-invertible elements, and there occurs a somewhat puzzling phenomenon: there are non-invertible elements $x$ such that $\det x\neq 0$. This fact is essentially due to the absolute convergence of the improper integral
\[\int_0^1\log s ds.\]

FK-determinant has been well-studied in many papers. In particular, the notion was extended to $C^*$-algebras in \cite{HS}. We refer the readers to \cite{Ha} for a recent survey, and make the following definition to proceed.\\

{\bf Definition.} In a unital $C^*$-algebra ${\mathcal B}$ with a normalized trace $tr$, an element $x$ will be called $tr$-singular if $\det x=0$.\\

This section will explicitly compute the FK-determinant for $R(z)=1+z_1a+z_2t$ (with respect to the representation $\lb_{D_{\infty}}$), and determine the points $z\in P(R)$ for which $R(z)$ is $tr$-singular. For convenience, we shall call these points $tr$-singular points in $P(R)$. For ${\mathcal B}$ as in the definition above, we consider an element $x\in GL({\mathcal B})$ that is in the path-connected component of the identity operator $I$ and a piecewise smooth path $x(s),\ 0\leq s\leq 1$, in $GL({\mathcal B})$ such that $x(0)=I$ and $x(1)=x$. In the case the integral
\[\int_0^1tr(x^{-1}(s)x'(s))ds\]
is independent of the path $x(s)$, the following quantity
\begin{equation*}
det_{tr}x:=\exp\big(\int_0^1tr(x^{-1}(s)x'(s))ds\big)
\end{equation*}
is well defined for $C^*$-algebras with trace (cf. \protect\cite{Ha,HS}).
Note that $det_{tr}$ can take on complex values, not just positive numbers. And it is shown in \cite{Ha} that
\begin{equation}
|det_{tr}x|=\det x.
\end{equation}
We observe that in fact (4.1) holds as long as $|det_{tr}(x)|$ (not $det_{tr}(x)$ ) is independent of the given path, or equivalently the integrals \[\int_0^1tr(x^{-1}(s)x'(s))ds\]
with respect to different paths differ by a purely imaginary number. This observation and Proposition 3.2 lead to the following theorem.

\begin{thm}
For $z=(z_1, z_2)\in \C^2$ and $R(z)=I+z_1a+z_2t$ with respect to the representation $\lb_{D_{\infty}}$, the Fuglede-Kadison determinant
\[\det R(z)=\exp\big(\frac{1}{4\pi}\int_0^{2\pi} \log |1-z_1^2-z_2^2-2z_1z_2\cos\theta|d\theta\big).\]
\end{thm}

\pf First, by Corollary 1.2 we know that $P^c(R)$ is path-connected. Let \[
\gamma=\{z(s): 0\leq s\leq 1\}\] be a closed piece-wise smooth path ($z(0)=z(1)$) in $P^c(R)$.
For simplicity, we let \[L_x(z)=1-z_1^2-z_2^2-2z_1z_2x.\] By Theorem 1.1 for the case $z_0=1$, for every fixed $-1\leq x\leq 1$ the function $L_x(z)$ does not vanish on $P^c(R)$. In particular, we have
\[L_x(z(s))=1-z_1^2(s)-z_2^2(s)-2z_1(s)z_2(s)x\neq 0,\ \forall s\in [0,\ 1],\] hence $L_x(\gamma)$ is a piecewise smooth path in the complex plane that does not run over $0$. We define the winding number $W(\gamma)$ of $\gamma$ around $P(R)$ as the winding number of $L_x(\gamma)$ around $0$, i.e.
\begin{align*}
W(\gamma)&=\frac{1}{2\pi i}\int_{L_x(\gamma)}\frac{1}{w}dw\\
&=\frac{1}{2\pi i}\int_{0}^1\frac{1}{L_x(z(s))}dL_x(z(s))\\
&=\frac{1}{2\pi i}\int_{0}^1\frac{d}{ds}\log L_x(z(s))ds.
\end{align*}
Observe that since $L_x(z)$ is linear in $x$, the above integrals show that $W(\gamma)$ is continuous with respect to $x$. But since $W(\gamma)$ is integer-valued, it is a constant with respect to the change of $x$, or in other words, the value of $W(\gamma)$ is independent of the choice of $x\in [-1,\ 1]$.
Using Proposition 3.2, Formula (3.7) and the above definition of winding number $W(\gamma)$, we can define the following coupling
\begin{align*}
\langle \gamma,\ tr\omega_R\rangle :&=\frac{1}{2\pi i}\int_{\gamma} tr \omega_R(z)\\
&=\frac{1}{2\pi i}\int_{\gamma}\partial \big(\frac{1}{4\pi}\int_0^{2\pi} \log (1-z_1^2-z_2^2-2z_1z_2\cos\theta)d\theta\big)\\
&=\frac{1}{4\pi}\int_0^{2\pi}\big(\frac{1}{2\pi i}\int^1_{0}\frac{d}{ds} \log (1-z_1^2(s)-z_2^2(s)-2z_1(s)z_2(s) \cos\theta)ds\big)d\theta\\
&=\frac{1}{4\pi}\int_0^{2\pi}\big(\frac{1}{2\pi i}\int_{0}^1\frac{d}{ds}\log L_{\cos \theta}(z(s))ds\big)d\theta\\
&=\frac{1}{4\pi}\int_0^{2\pi} W(\gamma)d\theta=\frac{W(\gamma)}{2}. \tag{4.2}
\end{align*}
Clearly, $(0,0)\in P^c(R)$. Since  $P^c(R)$ is path connected, for every $p\in  P^c(R)$ there is a piece-wise smooth path $\gamma (s)=(z_1(s), z_2(s)),\ 0\leq s\leq 1$, in $P^c(R)$ such that $\gamma (0)=(0, 0)$ and $\gamma(1)=z$, then again by Proposition 3.2 we have
\begin{align*}
\int_0^1tr \omega_R(\gamma(s))&=\int_0^1\partial \big(\frac{1}{4\pi}\int_0^{2\pi} \log (1-z_1(s)^2-z_2(s)^2-2z_1(s)z_2(s)\cos\theta)d\theta\big)\\
&=\frac{1}{4\pi}\int_0^{2\pi}\big(\int_0^1 \frac{d}{ds}\log (1-z_1(s)^2-z_2(s)^2-2z_1(s)z_2(s)\cos\theta)ds \big)d\theta\\
&=\frac{1}{4\pi}\int_0^{2\pi} \log (1-z_1^2-z_2^2-2z_1z_2\cos\theta)d\theta.
\end{align*}
If $\gamma'(s)$ is another such path connecting $(0,0)$ to $z$, then by (4.2),
\[\int_0^1tr \omega_R(\gamma(s))-\int_0^1 tr\omega_R(\gamma'(t))=2\pi i\langle \gamma-\gamma',\ tr\omega\rangle=W(\gamma-\gamma')\pi i,\]
where $\gamma-\gamma'$ stands for the closed path that goes from $0$ to $z$ along $\gamma$ and returns to $0$ along $\gamma'$. This indicates that the integrals
\[\int_0^1tr \omega_R(\gamma(s))=\int_0^1tr\big(R^{-1}(z(s))R'(z(s)\big)ds\]
with respect to different paths connecting $(0,0)$ to $z$ differ by a purely imaginary number, namely integer multiples of $\pi i$.
In conclusion, $|\exp \big(\int_0^1tr \omega_R(\gamma(s)\big)|$ is independent of the path, and by (4.1)
\begin{align*}
\det R(z)&=|\exp \big(\int_0^1tr \omega_R(\gamma(s)\big)|\\
&=\exp\big(Re\frac{1}{4\pi}\int_0^{2\pi} \log (1-z_1^2-z_2^2-2z_1z_2\cos\theta)d\theta\big)\\
&=\exp\big(\frac{1}{4\pi}\int_0^{2\pi} \log |1-z_1^2-z_2^2-2z_1z_2\cos\theta|d\theta\big).\tag{4.3}
\end{align*}
Since for each fixed $z\in \C^2$ and all $\theta\in [0,\ 2\pi]$ \[\log (1+(|z_1|+|z_2|)^2)\geq\log |1-z_1^2-z_2^2-2z_1z_2\cos\theta|\geq -\infty,\] and the integral in (4.3) is convergent or equal to $-\infty$,
formula (4.3) extends to all $z\in \C^2$, and the theorem is established.
\zb

There is an interesting special case of Theorem 4.1. Consider the quadratic surface $S\subset \C^2$ defined by $1-z_1^2-z_2^2=0$. Then on $S$,
\[1-z_1^2-z_2^2-2z_1z_2\cos \theta=-2z_1z_2\cos \theta,\]
which vanishes at $\theta=\pi/2$ and $3\pi/2$. So by Theorem 1.1, $S\subset P(R)$. To proceed, we first
verify the formula
\begin{equation*}
\int_0^{\pi/2} \log \cos\theta d\theta=-\frac{\pi}{2}\log 2.\tag{4.4}
\end{equation*}
Denoting the integral by $M$ and using the fact $\cos \theta=\sin (\frac{\pi}{2}-\theta)$, one checks that
\begin{align*}
2M&=\int_0^{\pi/2} \log \cos\theta d\theta+\int_0^{\pi/2} \log \sin (\frac{\pi}{2}-\theta) d\theta\\
&=\int_0^{\pi/2} \log \cos\theta d\theta+\int_0^{\pi/2} \log \sin\theta d\theta\\
&=\int_0^{\pi/2} \log\sin(2\theta)-\log 2d\theta\\
&=\frac{1}{2}\int_0^{\pi} \log\sin(\theta)d\theta-\frac{\pi}{2}\log 2\\
&=\frac{1}{2}\left(\int_0^{\pi/2} \log\sin(\theta)d\theta+\int_0^{\pi/2} \log\sin(\theta+\frac{\pi}{2})d\theta\right)-\frac{\pi}{2}\log 2.
\end{align*}
Formula (4.4) then follows easily from the fact $\sin (\theta+\frac{\pi}{2})=\cos\theta$.

On $S$, using the formula (4.4) we compute that
\begin{align*}
\det R(z)&=\exp\big(\frac{1}{4\pi}\int_0^{2\pi} \log |2z_1z_2|+\log|\cos\theta| d\theta\big)\\
&=\sqrt{2|z_1z_2|}\exp \big(\frac{1}{4\pi}\int_0^{2\pi} \log |\cos\theta| d\theta\big)\\
&=\sqrt{2|z_1z_2|}\exp \big(\frac{1}{\pi}\int_0^{\pi/2} \log \cos\theta d\theta\big)\\
&=\sqrt{2|z_1z_2|}/ \sqrt{2}\\
&=\sqrt{|z_1z_2|}.\tag{4.5}
\end{align*}
This confirms that $\det R(z)$ may be nonzero even though $R(z)$ is not invertible.

The $tr$-singular points in $P(R)$ (points for which $\det R(z)=0$) are now easy to determine. First, if $z_1z_2\neq 0$, then after factoring out $2z_1z_2$, the integral in (4.3) is of the form
\[\int_0^{2\pi} \log |\beta-\cos\theta|d\theta,\]
which is well-known to be convergent. Hence $\det R(z)\neq 0$ in this case. If $z_1z_2=0$, then the integral in (4.3) diverges only if $1-z_1^2-z_2^2=0$. Solving these two equations, we have
\begin{corr}
For $R(z)=1+z_1a+z_2t$, the set of $tr$-singular points in $P(R)$ is $\{(\pm1,\ 0),\ (0,\ \pm1)\}$.
\end{corr}

 We conclude this section with a few observations.

{\bf 1}. Connection between FK-determinant and Mahler measure has been noted in quite a few recent papers (for example \cite{De,Ha,Li,Sc}). For references on Mahler measure, we refer readers to \cite{EW}. This connection is clean and explicit in the case here. For a complex polynomial
\[P(w)=a(w-\alpha_1)(w-\alpha_2)\cdots (w-\alpha_n),\]
its Mahler measure is defined as
\[M(P)=|a|\prod_{|\alpha_j|\geq 1}|\alpha_j|,\]
and it is well-known that by Jensen's formula (\cite{Ru})
\[M(P)=\exp\big(\frac{1}{2\pi}\int_0^{2\pi}\log|P(e^{i\theta})|d\theta\big).\]
If we let $P_z(w)=w(1-z_1^2-z_2^2)-z_1z_2(w^2+1)$, then on the unit circle $w=e^{i\theta}$ and
\[|P_z(w)|=|(1-z_1^2-z_2^2)-z_1z_2(w+\overline{w})|=|1-z_1^2-z_2^2-2z_1z_2\cos\theta|,\]
and Theorem 4.1 has the following immediate consequence.
\begin{corr}
$\det R(z)=\sqrt{M(P_z)}.$
\end{corr}
$M(P_z)$ can of course be computed in this case. But the expression is not as clean.
Observe that the set $\{(\pm1,\ 0),\ (0,\ \pm1)\}$ is precisely the set on which the polynomial $P_z(w)$ in Corollary 4.3 is constant $0$.

{\bf 2}. The coupling in (4.2) yields some information about the singular homology group $H_1(P^c(R),\mathbb Z)$. Note that since $P^c(R)$ is connected, by Hurewicz theorem $H_1(P^c(R),\mathbb Z)$ is isomorphic to the abelianization of the homotopy group $\pi_1(P^c(R))$.

\begin{corr}
The coupling with $tr\omega_R$ in (4.2) defines an isomorphism from $H_1(P^c(R), {\mathbb Z})$ onto $\frac{1}{2}{\mathbb Z}$.
\end{corr}
\pf First, differentiating both sides of $R^{-1}(z)R(z)=I$, we have
\[\left(dR^{-1}(z)\right)R(z)+R^{-1}(z)dR(z)=0,\]
and hence
\[dR^{-1}(z)=-R^{-1}(z)(dR(z))R^{-1}(z).\]
Consider \[\omega_R(z)=R^{-1}(z)dR(z)=(I+z_1a+z_2t)^{-1}(adz_1+tdz_2),\ \ z\in P^c(R).\] Then by above computation one sees that
\begin{align*}
d\omega_R(z)&=(dR^{-1}(z))\wedge dR(z)\\
&=-\omega_R(z)\wedge \omega_R(z)\\
&=-\left(R^{-1}(z)aR^{-1}(z)t-R^{-1}(z)tR^{-1}(z)a\right)dz_1\wedge dz_2.
\end{align*}
Hence
\begin{align*}
d[tr\omega_R(z)]&=tr[d\omega_R(z)]\\
&=-tr \left(R^{-1}(z)aR^{-1}(z)t-R^{-1}(z)tR^{-1}(z)a\right)dz_1\wedge dz_2=0,
\end{align*}
meaning that $tr\omega_R$ is a closed holomorphic $1$-form. If $\tilde{\gamma}$ is a closed piece-wise smooth path that is homological to $\gamma$, i.e. $\gamma-\tilde{\gamma}=\partial \Delta$
for some singular $2$-simplex $\Delta$ (where $\partial$ stands for the boundary map here), then an argument using Stokes theorem shows that
\begin{align*}
\langle \gamma,\ tr\omega_R\rangle -\langle \tilde{\gamma},\ tr\omega_R\rangle&=\int_{\gamma}tr\omega-\int_{\tilde{\gamma}}tr\omega_R\\
&=\int_{\Delta} d[tr\omega_R]=0.
\end{align*}
This means that the coupling with $tr\omega_R$ in (4.2) gives rise to a homomorphism $\kappa$ from the singular homology group $H_1(P^c(R), {\mathbb Z})$ into $\frac{1}{2}{\mathbb Z}$ defined by
\[\kappa([\gamma])=\langle \gamma,\ tr\omega_R\rangle.\]
Further, if $\langle \gamma,\ tr\omega_R\rangle =0$, then the winding number $W(\gamma)=0$ by (4.2), indicating that $\gamma$ does not wind around $P(R)$ and hence is homological in $P^c(R)$ to a point, i.e. $[\gamma]=0$ in $H_1(P^c(R), {\mathbb Z})$. This shows that $\kappa$ is injective.

To check the range of the map $\kappa$, we consider the path
\[\gamma=\{z(s)=(1+e^{2\pi is}/2,\ 0),\ 0\leq s\leq 1\}.\]
Then for every $x\in [-1,\ 1]$, \[L_x(z(s))=1-(1+e^{2\pi is}/2)^2=-e^{2\pi is}-e^{4\pi is}/4\neq 0\]
for all $s\in [0,\ 1]$. Hence by Theorem 1.1 for the case $z_0=0$, the path $\gamma$ is inside $P^c(R)$.
By the definition in the proof of Theorem 4.1, $\gamma$'s winding number
\begin{align*}
W(\gamma)&=\frac{1}{2\pi i}\int_{0}^1\frac{d}{ds}\log L_x(z(s))ds\\
&=\frac{1}{2\pi i}\int_{0}^1\frac{-2(1+e^{2\pi is}/2)\pi i e^{2\pi is}ds}{-e^{i2\pi s}-e^{4\pi is}/4}\\
&=\int_{0}^1\frac{1+e^{2\pi is}/2ds}{1+e^{2\pi is}/4}.
\end{align*}
Using the geometric series for $(1+e^{2\pi is}/4)^{-1}$ and direct computation of integrals, we have $W(\gamma)=1$. Then by (4.2), we have
\[\kappa([\gamma])=\langle \gamma,\ tr\omega_R\rangle=\frac{W(\gamma)}{2}=\frac{1}{2}.\]
Therefore,
\[\kappa(n[\gamma])=\frac{nW(\gamma)}{2}=n/2,\ \ n\in {\mathbb Z}.\]
This completes the proof.\zb

Although $\frac{1}{2}{\mathbb Z}$ is isomorphic to ${\mathbb Z}$, the factor $1/2$ here tells more about the structure of the group $D_{\infty}$.

{\bf 3}. Note also that the classical spectra $\sigma(a)=\sigma(t)=\{\pm 1\}$. It suggests that $tr$-singular points may have a direct connection with the classical spectra.

{\bf 4}. For $R(z)=1+z_1a+z_2t$, by \cite{DY} the $(1,1)$-form \[\frac{-i}{2}tr\big(\omega^*_R(z)\wedge \omega_R(z)\big)\] induces a natural Riemannian metric on $P^c(R)$. The metric has singularities at points in $P(R)$. It follows from a result in \cite{DY} that the $tr$-singular points in $P(R)$, namely $\{(\pm 1, 0),\ (0,\pm 1)\}$, are outside of the completion of $P^c(R)$ under this metric. It is not clear if the completion of $P^c(R)$ is equal to $\C^2\setminus \{(\pm 1, 0),\ (0,\pm 1)\}$.

{\bf 5}. Likewise, if we consider $R(z)=z_0+z_1a+z_2t$, then it is not hard to check that
\[\det R(z)=\exp\big(\frac{1}{4\pi}\int_0^{2\pi} \log |z_0^2-z_1^2-z_2^2-2z_1z_2\cos\theta|d\theta\big).\]

\section{On the finite dihedral group}

For the finite dihedral group $D_{n}=<a,\ t\ |\ a^2=t^2=(at)^n=1>$, where $n\in {\mathbb N}$, $H=< at>$ is
a cyclic normal subgroup of order $n$. The GNS construction for its group $C^*$-algebra $C^*(H)$
gives the finite dimensional Hilbert space \[L_n=\{f= \sum_{0}^{n-1} f_j(at)^j\},\]
with inner product
\[\langle g,\ f\rangle=tr(f^*g)=\sum_{0}^{n}\overline{f_j}g_j.\]
Similar to (1.1), the left regular representation for $a$ and $t$ on $L_n\oplus L_n$ is given by
\begin{eqnarray}
\lb_{D_n}(a)=\left(
      \begin{array}{cc}
        0 & T  \\
        T^{*} & 0 \\
      \end{array}
    \right),\quad
\lb_{D_n}(t)=\left(
      \begin{array}{cc}
        0 & 1  \\
        1 & 0 \\
      \end{array}
    \right),
\end{eqnarray}
where $T$, with respect to the orthonormal basis $\{1,\ at,\ \cdots,\ (at)^{n-1}\}$, is the unitary
\[\begin{bmatrix}
        0 & 0 & 0 \cdots & \cdots & 1 \\
        1 & 0 & 0 \cdots & \cdots & 0 \\
        0 & 1 & 0 \cdots & \cdots & 0 \\
        \vdots & \vdots & \vdots & \vdots & \vdots\\
        0 & 0 & 0 \cdots & 1 & 0
      \end{bmatrix}_{n\times n}.\]
Clearly $T^n=1$ and its spectrum $\sigma(T)=\{e^{\frac{2k\pi i}{n}},\ 0\leq k\leq n-1\}$. Denoting $\frac{2k\pi }{n}$ by $\theta_k$, $T$'s projection-valued spectral measure $E(e^{i\theta})$ is equal to the projection onto the sum of eigenspaces corresponding to the eigenvalues $\theta_k\leq \theta $. Since $trI=1$,
\[trE(e^{i\theta})=\frac{1}{n}|\{k:\ \theta_k\leq \theta\}|.\]
 By arguments similar to that in Sections 3 and 4, we have
\begin{corr}
Consider $R(z)=z_0+z_1a+z_2t$. Then with respect to the representation of $D_n$ given by $\lb_{D_n}$, we have
\begin{align*}
 (a)&\ \ P(R)=\bigcup_{k=0}^{n-1}\{z\in \C^3:\ z_0^2-z_1^2-z_2^2-2z_1z_2\cos \theta_k=0\}.\\
 (b)&\ \ trR^{-1}(z)=\frac{1}{n}\sum_0^{n-1}\frac{z_0}{z_0^2-z_1^2-z_2^2-2z_1z_2\cos \theta_k}.\\
 (c)&\ \ tr\omega_R(z)=\partial \big(\frac{1}{2n}\sum_{k=0}^{n-1} \log (z_0^2-z_1^2-z_2^2-2z_1z_2\cos\theta_k)\big).
\end{align*}
\end{corr}

By the matrix expression of $T$ and the representation (5.1), $R(z)=z_0+z_1a+z_2t$ is a $2n\times 2n$ matrix. Its matrix determinant, denoted by  $det R(z)$, can be computed, but the process is complicated if computed directly. Instead, we can use the idea leading to Theorem 4.1. By Corollary 5.1(c)
\begin{align*}
\det R(z)&=\exp\big(\frac{1}{2n}\sum_{k=0}^{n-1} \log |z_0^2-z_1^2-z_2^2-2z_1z_2\cos\theta_k|\big)\\
&=\big(\prod_{k=0}^{n-1} |z_0^2-z_1^2-z_2^2-2z_1z_2\cos\theta_k|\big)^{1/2n}.
\end{align*}
For a $k\times k$ matrix algebra, it is known that the Fuglede-Kadison determinant and the regular determinant are related by the formula $\det A=|det A|^{1/k}$. In here, since $R(z)$ is $2n\times 2n$, we have
\begin{align*}
|det R(z)|=(\det R(z))^{2n}=\prod_{k=0}^{n-1}|z_0^2-z_1^2-z_2^2-2z_1z_2\cos \frac{2k\pi}{n}|.
\end{align*}
Since $detR(1,0,0)=det(I)=1$, the next corollary follows.

\begin{corr}
With respect to the representation $\lb_{D_n}$,
\[det R(z)=\prod_{k=0}^{n-1}(z_0^2-z_1^2-z_2^2-2z_1z_2\cos \frac{2k\pi}{n}).\]
\end{corr}
This clearly verifies Corollary 5.1(a).\\

{\bf Remark.} We observe that since $D_n$ is a finite group, it is amenable, and hence by the Hulanicki-Reiter theorem in Section 2, the representation $\lb_{D_n}$ in (5.1) (which is unitarily equivalent to the left regular representation) is maximal in the sense of weak contaiment. The next corollary then follows from Proposition 2.2 and Corollary 5.1(a).

\begin{corr}
If $A$ and $T$ are elements in a $C^*$-algebra $B$ such that $A^2=T^2=(AT)^n=I$, and $R_B(z)=z_0I+z_1A+z_2T$, then $P(R_B)$ is a subset in
\[\bigcup_{k=0}^{n-1}\{z\in \C^3:\ z_0^2-z_1^2-z_2^2-2z_1z_2\cos \theta_k=0\}.\]
\end{corr}

\begin{exam}
In Pedersen's representation in Section 2, if we let $x=1/2$, then
\[a=\begin{pmatrix} 0 & -1 \\ -1 & 0 \end{pmatrix},\
  t=\begin{pmatrix} -1 & 0 \\ 0 & 1 \end{pmatrix},\]
are non-commuting self-adjoint unitaries such that $(at)^4=1$. So $a$ and $t$ generate $D_4$. But as $2\times 2$ matrices, $R(z)=z_0+z_1a+z_2t$ is not invertible if and only if its ordinary determinant
\[det R(z)=z_0^2-z_1^2-z_2^2=0.\]
So in this case $P(R)$ consists of a single surface $\{z_0^2-z_1^2-z_2^2=0\}$, while with respect to the representation of $D_4$ given by (5.1), $P(R)$ has $4$ pieces by Corollary 5.1(a) (the piece $\{z_0^2-z_1^2-z_2^2=0\}$ has multiplicity $2$.).
\end{exam}

\section{Self-similarity}

In Section 2, we pointed out that for weakly equivalent representations $\pi$ and $\rho$ we have $P(A_{\pi})= P(A_{\rho})$. It is an appealing question that in which cases the converse is true. We will address this issue in these two sections for the case $G=D_{\infty}$. We are going to show that $C^*(D_{\infty})$ can be realized as a self-similar algebra in the sense of definition 2.1 in \cite{GN07}. This will be done via computation of the joint spectrum of two involutive operators associated with a certain self-similar realization of $D_{\infty}$ and, surprisingly, will be linked to the $3$-generated $2$-group $\G$ of intermediate growth constructed by the first author in \cite{Gr80} (cf. (0.1)) and studied in \cite{Gr84,Gr89,Gr05} and many other articles.

This section gives a brief introduction to self-similarity. The concept of self-similarity entered algebra via group theory, namely via a careful study of self-similarity properties of $\G$ and its relatives.

Let $X=\{x_1,\ x_2,\ \cdots,\ x_d\}$ be a finite alphabet of cardinality $d\geq 2$, $X^*$ be the set of all finite words over $X$. Then
\begin{equation}
X^*=\bigcup_{n=0}^{\infty}X^n,
\end{equation}
 where $X^n$ is the set of words of length $n$. A self-similar action of a group $G$ is an action of $G$ on $X^*$ such that for every $g\in G$ and every $x\in X$ there exist $h\in G$ and $y\in X$ such that
\begin{equation}
g(xw)=yh(w)
\end{equation}
for all $w\in X^*$. Observe that such an action preserves the length of the words, so each $X^n$ is $G$-invariant.

Let $T=T(X)$ be a $d$-regular rooted tree with set $V$ of vertices $X^*$ and the set $E=\{(w,\ wx)|\ w\in X^*,\ x\in X\}$ of edges. For a vertex $v\in V$, $T_v$ is the subtree of $T$ with root at $v$. Figure (7.1) shows the first $3$ levels of an infinite binary rooted tree for $X=\{0,\ 1\}$. For example, $T_{01}$ is the subtree with root at $01$.

It is easy to see that the self-similar action $(G,\ X^*)$ induces an action of $G$ on $T(X)$ by automorphisms. Then the root (corresponding to empty word) and the levels $V_n=X^n$ of the tree are $G$-invariant, and the maximum of transitivity that may occur is the level transitivity (i.e. the transitivity of action on each level). For each fixed $n\in {\mathbb N}$ and $g\in G$ one has a decomposition $g=(g_1,\ g_2,\ \cdots,\ g_{d^n})\sigma$, where $\sigma\in Sym(d^n)$ and $g_i\in Aut(T_i)$. Here $T_i,\ 1\leq i\leq d^n$, is the subtree of $T$ with root at the $i$-th vertex of level $n$ (see Figure (7.1)). Elements $g_1,\ g_2,\ \cdots,\ g_{d^n}$ are called sections of $g$ in corresponding vertices and the action $(G,\ X^*)$ is self-similar if and only if for every $n\in {\mathbb N},\ 1\leq i\leq d^n$ the section $g_i$ belongs to $G$ after canonical identification of $T_i$ with $T$. It is not hard to see that  $(G,\ X^*)$ is self-similar if and only if sections $g_1,\ g_2,\ \cdots,\ g_{d}$ in vertices of the first level belong to $G$ (after identification $T_i\cong T,\ i=1,\ 2,\ \cdots,\ d$).

A group $G$ is called self-similar if for some $X$ it has a faithful self-similar action on $X^*$. An alternative definition of self-similarity is that $G$ is isomorphic to a group $G({\mathcal A})$ generated by non-initial Mealy type automaton ${\mathcal A}$ over finite alphabet. For more details, see \cite{GNS} or \cite{Ne05}.

Every self-similar group is residually finite (i.e. approximated by finite groups) but not vice versa. We call $G$ an $(m,\ d)$-group if $G\cong G({\mathcal A})$ where ${\mathcal A}$ is an automaton with $m$ states over an alphabet with $d$ symbols. There are six $(2,\ 2)$ groups listed in \cite{GNS} and $D_{\infty}$ is in this list. Two essentially different realizations of $D_{\infty}$ as $(2,\ 2)$-group are:
 
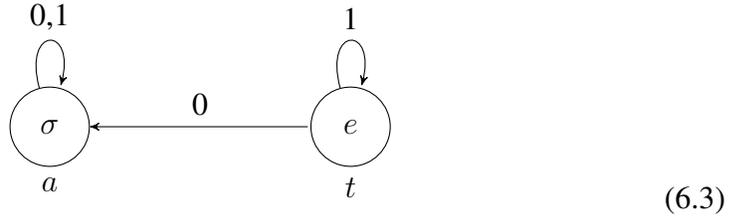
\begin{align}
\begin{tikzpicture}[>=stealth',shorten >=1pt,auto,node distance=4cm]
 \node[state,label=below:$a$] (S) {$\sigma$};
\node[state,label=below:$t$]         (e) [right of=S] {$e$};
\path[<-] (S)  edge [loop above] node {0,1} (S)
            edge              node {0} (e)
        (e) edge [loop above] node {1} (e);
\end{tikzpicture}
\end{align}
where $a=(a,\ a)\sigma,\ t=(a,\ t)$, and
 \begin{align}
& 
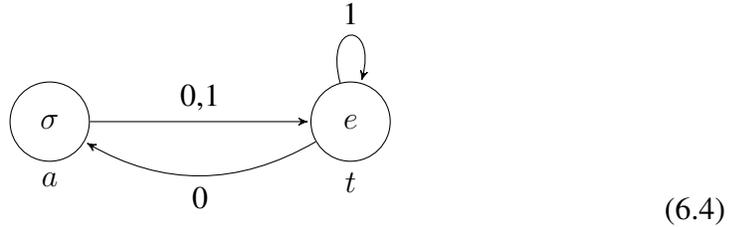
\begin{tikzpicture}[>=stealth',shorten >=1pt,auto,node distance=4cm]
\node[state,label=below:$a$] (S) {$\sigma$};
\node[state,label=below:$t$]         (e) [right of=S] {$e$};
\path[->] (S)
            edge              node {0,1} (e)
        (e) edge [loop above] node {1} (e);
  \path[->]  (e) edge [bend left]  node {0} (S);
\end{tikzpicture}
\end{align}
where $a=(t,\ t)\sigma,\ t=(a,\ t)$. Here $e$ is the identity map on $\{0,\ 1\}$ and $\sigma$ is the involution.

There is no complete classification of $(3,\ 2)$-groups yet, although much information is collected in \cite{BG..S}, and it is known that up to isomorphism there is no more than 115 such groups. The following realization of $D_{\infty}$ as a $(3,\ 2)$-group is used in the sequel. $D_{\infty}$ is isomorphic to $G({\mathcal A})$ where ${\mathcal A}$ is the automaton

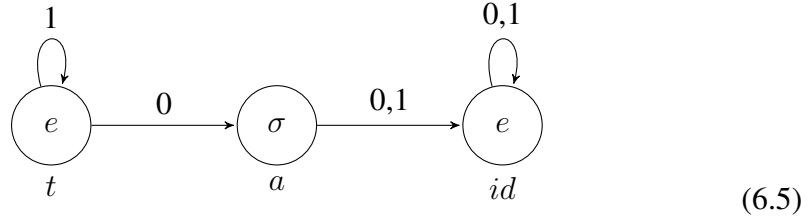
\begin{align}
\begin{tikzpicture}[>=stealth',shorten >=1pt,auto,node distance=3cm]
 \node[state,label=below:$t$] (S) {$e$};
\node[state,label=below:$a$]         (T) [right of=S] {$\sigma$};
\node[state,label=below:$id$]         (U) [right of=T] {$e$};
\path[->] (S)  edge [loop above] node {1} (S);
 \path[->]            (S) edge              node {0} (T);
\path[->] (T)       edge              node {0,1} (U);
\path[->] (U) edge [loop above] node {0,1} (U);
\end{tikzpicture}
\end{align}
or $D_{\infty}=<a,\ t>$ where $a$ and $t$ are automorphisms of binary tree satisfying the recursive relation
\begin{align}
a=\sigma,\ \ \ \ t=(a,\ t),
\end{align}
where $\sigma$ is the involution. This realization corresponds to the automaton 2874 in \cite{BG..S}.

Now we are going to define self-similar representations. Let $H$ be a separable Hilbert space. By a $d$-similarity of $H$ we mean (following \cite{GN07}) an isomorphism of Hilbert spaces $\psi:\ H\longrightarrow H^d$. Fix such a $\psi$ and define $T_1,\ T_2,\ \cdots,\ T_d\in B(H)$ (i.e. bounded linear operators acting on $H$) by
\begin{equation}
T_k(\xi)=\psi^{-1}(0,\ \cdots, \ 0,\ \xi,\ 0,\ \cdots,\ 0),
\end{equation}
where $\xi\in H$ and is at the $k$-th coordinate on the right-hand side. Let $G$ act self-similarly on $T(X)$ ($|X|=d$). A unitary representation $\rho$ of $G$ on $H$ is said to be self-similar (with respect to the $d$-similarity $\psi$) if \[\rho(g)T_x=T_y\rho(h)\]
whenever $g(xw)=yh(w)$ for all $w\in X^*$ (i.e. whenever (6.1) holds).

For a self-similar representation $\rho:\ G\longrightarrow B(H)$, every operator $\rho(g),\ g\in G$ can be written, with respect to the $d$-similarity $\psi$, as a $d\times d$ matrix $\rho(g)=(A_{yx})_{x,y\in X}$, where
\[A_{yx}=\left\{\begin{array}{ll}
                     \rho(g|_x) & \ \text{if}\ g(x)=y\\
                     0 & \text{otherwise,}
                    \end{array}
                \right. \]
and $g|_x$ denotes the section of $g$ at vertex $x$ of the first level of $T(X)$.

Extension of $\rho$ to group algebra $\C[G]$ leads to a self-similar representation of $\C[G]$, and induces a map $\phi :\ \C[G]\longrightarrow M_d(\C[G])$. Here $M_d(\textit{R})$ denotes the ring of matrices over ring $\textit{R}$. The map $\phi$ is injective on $G$ (but not on $\C[G]$, cf. \cite{GN07}), and for $g\in G$
\begin{equation}
\phi(g)=(A_{yx}).
\end{equation}
Moreover, it extends to a homomorphism
\begin{equation}
C^*_{\rho}(G)\longrightarrow M_d(C^*_{\rho}(G)),
\end{equation}
which we also denote by $\phi$. This $\phi$ is injective, since it implements the equivalence of the representation $\rho$ with the representation $\psi \circ \rho \circ \psi^{-1}$.

Following \cite{GN07} we call a completion of $\C[G]$ self-similar if it is the completion with respect to the norm defined by self-similar representation, and call the corresponding $C^*$-algebra $C^*_{\rho}(G)$ self-similar.

\section{Self-similar realization of $C^*(D_{\infty})$}

There is a canonical self-similar representation associated with a self-similar group $G$ acting (self-similarly) on $T(X)$. Let $\partial T=X^{\mathbb N}$ be the boundary of the tree $T=T(X)$ (supplied with Tychonoff topology). Let $\mu=\{1/d,\ 1/d,\ \cdots,\ 1/d\}^{\otimes {\mathbb N}}$ be a uniform Bernoulli measure on $\partial T$ and let $H=L^2(\partial T,\ \mu)$. The measure $\mu$ is invariant with respect to arbitray automorphisms of $T$, in particular, $\mu$ is $G$-invariant. The Koopman representation $\rho$ is defined by
\[\big(\rho(g)f\big)(x)=f(g^{-1}x)\]
for $g\in G$, $f\in H$. It is well-known that $\rho$ is a self-similar representation (cf. \cite{BG,GN07}).
Now we are in position to state the next theorem.

\begin{thm}
Let $D_{\infty}$ be acting self-similarly on a rooted tree according to the recursion relation (6.6). Let $\rho$ be the corresponding Koopman representation. Then $C^*_{\rho}(D_{\infty})$ is isomorphic to $C^*(D_{\infty})$, and hence $C^*_{\rho}(D_{\infty})$ is a self-similar realization of $C^*(D_{\infty})$.
\end{thm}

This result follows from the following

\begin{thm}
Let $\lb_{D_{\infty}}$ and $\rho$ be the left regular representation and the Koopman representation of $D_{\infty}$, respectively. Then

a) $\lb_{D_{\infty}}$ and $\rho$ are weakly equivalent;

b) the projective joint spectrum of the triple $(\rho(1),\ \rho(a),\ \rho(t))$ is the same as that in the case of regular representation, i.e. it equals the set
\[\bigcup_{-1\leq x\leq 1}\{z\in \C^3:\ z_0^2-z_1^2-z_2^2-2z_1z_2x=0\}.\]
\end{thm}
Here $a$ and $t$ are generators of $D_{\infty}$ corresponding to states of the automaton given by Figure (6.5).

\pf We begin with the proof of part b). Figure (7.1) shows the first three levels of an infinite binary tree. Let $T_i,\ i=0,\ 1$ be the subtrees of the rooted binary tree with roots at vertices of the first level, and $\mu_i,\ i=0,\ 1$ be the normalized restrictions of $\mu$ on $\partial T_0,\ \partial T_1$, and $H_i=L^2(\partial T_i,\ \mu_i)$, $H=H_0\oplus H_1$.
\begin{align}
\begin{tikzpicture}[->,>=stealth',level/.style={sibling distance = 5cm/#1,
  level distance = 1.5cm}]
\node [arn_r] {}
    child{ node [arn_r] {$0$}
            child{ node [arn_r] {$00$}}
            child{ node [arn_r] {$01$}}
           }
    child{ node [arn_r] {$1$}
            child{ node [arn_r] {$10$}}
            child{ node [arn_r] {$11$}}
	   }
;
\end{tikzpicture}
\end{align}

There are canonical identifications of $T_0$ and $T_1$ with $T$ which induce identification $H_i\cong H$ and hence $2$-similarity $\psi:\ H=H_0\oplus H_1\longrightarrow H\oplus H.$ The relations (6.8) and the embedding (6.9) in our case lead to the relations
\[
\phi(\rho(a))=\left(\begin{matrix}
                      0 & I\\
                      I & 0
                     \end{matrix}\right),\ \ \ \ \phi(\rho(t))=\left(\begin{matrix}
                                                                      \rho(a) & 0\\
                                                                      0 & \rho(t)
                                                                     \end{matrix}\right),\]
which we will simply write as
\begin{equation}
 \rho(a)=\left(\begin{matrix}
                      0 & I\\
                      I & 0
                     \end{matrix}\right),\ \ \ \ \rho(t)=\left(\begin{matrix}
                                                                      \rho(a) & 0\\
                                                                      0 & \rho(t)
                                                                     \end{matrix}\right).
\end{equation}
Here $I$ is the identity operator. Now consider $R_{\rho}(z)=\rho(1)+z_1\rho(a)+z_2\rho(t)$, where $z\in \C^2$, and we shall compute the joint spectrum
\[P(R_{\rho})=\{z\in \C^2:\ R_{\rho}(z) \ \text{is not invertible in } C^*_{\rho}(D_{\infty})\}.\]
The method is very different from that in the case of the regular representation. Indeed the result follows from computations in \cite{BG} on the group $\G$ in (0.1).

Let $\G=<a,\ b,\ c,\ d>$ be a self-similar group acting on a binary rooted tree given by recursive relations
\begin{equation}
a=\sigma,\ b=(a,\ c),\ c=(a,\ d),\ d=(1,\ b),
\end{equation}
or equivalently the group generated by the $5$-state automaton

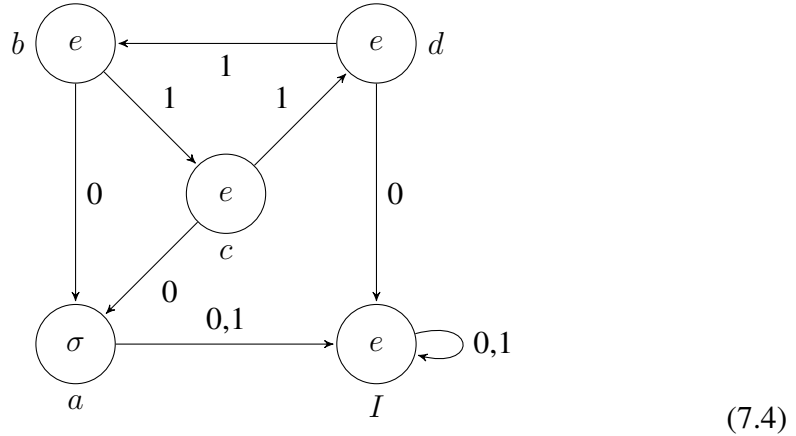
\begin{equation}
\begin{tikzpicture}[>=stealth',shorten >=1pt,auto,node distance=3cm]
  \node[state,label=left:$b$] (A) at ( -2,2) [shape=circle,draw]{$e$};
  \node[state,label=right:$d$] (B) at ( 2,2) [shape=circle,draw]{$e$};
  \node[state,label=below:$c$] (C) at ( 0,0) [shape=circle,draw]{$e$};
  \node[state,label=below:$a$] (D) at ( -2,-2) [shape=circle,draw] {$\sigma$};
  \node[state,label=below:$I$] (E) at (2,-2) [shape=circle,draw]{$e$};
\path[->]
        (B)   edge                    node {1}            (A)
        (A)   edge                    node {0}            (D)
        (A)   edge                    node {1}            (C)
        (C)   edge                    node {0}            (D)
        (C)   edge                    node {1}            (B)
        (D)   edge                    node {0,1}         (E)
        (B)   edge                    node {0}            (E)
        (E)   edge  [loop right]   node {0,1}            (E);
\end{tikzpicture}
\end{equation}

Observe that $a\in D_{\infty}\cap \G$. $\G$ can be described by generators and relations as in (0.1).
It is well-known that $\G$ is an infinite torsion $2$-group (group of Burnside type), has intermediate growth (between polynomial and exponential) and has many other unusual properties.

Now let $\pi:\ \G\longrightarrow U(H)$ be the Koopman unitary representation, and $H=L^2(\partial T,\ \mu)$ as before. Consider the pencil $Q(\lb,\mu)=\pi(-\lb a +b+c+d-(\mu+1)1)$, and let \[P(Q)=\{(\lb,\mu)\in \C^2:\ Q(\lb,\mu)\ \text{is not invertible in}\ C^*_{\pi}(\G)\}.\] The joint real spectrum $P(Q)\cap {\mathbb R}^2$ was computed in \cite{BG}, and the results extend naturally to $\C^2$. We are going to relate
 the spectra of $R_{\rho}(z)$ and  $P(Q)$. Indeed, if we let $u=(b+c+d-1)/2$ and rescale $Q$ we have
\[\frac{Q}{-\mu}=\frac{\lb}{\mu}\pi(a)+\frac{-2}{\mu}\pi(u)+1.\]
Then letting $z_1=\frac{\lb}{\mu}$ and $z_2=\frac{-2}{\mu}$, we have an expression similar to $R_{\rho}(z)$. We shall show that indeed $u$ can be identified with $t$. By (0.1) (or (7.3)), one sees that
\[bc=d=d^{-1}=cb,\ \ cd=b=b^{-1}=dc, \ \ bd=c=c^{-1}=db,\]
which indicates that elements $1,\ b,\ c,\ d$ constitute the Klein group ${\mathbb Z}_2\oplus {\mathbb Z}_2$, and in particular the element $u \in \C[\G]$ satisfies the relation
\[u^2=\frac{4+2bc+2cd+2bd-2b-2c-2d}{4}=1.\]
Therefore the group algebra $\C[\G]$ contains the dihedral group $\overline{D}:=< a,\ u>$ via the isomorphism defined by $a\to a$, $u\to t$. We shall check that $au$ is of infinite order and hence $\overline{D}$ is the infinite dihedral group.

\begin{lem}
The operators $\rho(t)$ and $\pi(u)$ coincide.
\end{lem}
\pf The self-similarity relations \[a=\sigma,\ b=(a,\ c),\ c=(a,\ d),\ d=(1,\ b)\]
lead to the following operator recursions
\begin{align*}
&\pi(a)=\left(\begin{matrix}
                      0 & I\\
                      I & 0
                     \end{matrix}\right),\ \ \ \ \pi(b)=\left(\begin{matrix}
                                                                      \pi(a) & 0\\
                                                                      0 & \pi(c)
                                                                     \end{matrix}\right)\\
&\pi(c)=\left(\begin{matrix}
                      \pi(a) & 0\\
                      0 & \pi(d)
                     \end{matrix}\right),\ \ \ \ \pi(d)=\left(\begin{matrix}
                                                                      I & 0\\
                                                                      0 & \pi(b)
                                                                     \end{matrix}\right).
\end{align*}
Therefore
\[\pi(u)=\pi\big(\frac{b+c+d-1}{2}\big)=\left(\begin{matrix}
                                                                      \pi(a) & 0\\
                                                                      0 & \pi(u)
                                                                     \end{matrix}\right),\]
and we see that $\pi(u)$ satisfies the same operator recursion as
\begin{equation}
\rho(t)=\left(\begin{matrix}
                     \rho(a) & 0\\
                             0 & \rho(t)
                    \end{matrix}\right)
\end{equation}
and
\begin{equation}
\rho(a)=\pi(a)=\left(\begin{matrix}
                                                                      0 & I\\
                                                                      I & 0
                                         \end{matrix}\right).
\end{equation}

It is natural to expect that $\pi(u)=\rho(t)$ as operators on $H=L^2(\partial T,\ \mu)$. This is indeed the case, but additional arguments are needed to justify this claim. Let $\pi_n$ be the permutational representation of ${\overline{D}}$ in $l^2(V_n)$, where as before $V_n$ is the $n$-th level of the tree $T$. The space $l^2(V_n)$ can be naturally identified with the subspace
\[H_n=span\{\chi_{E_j^{(n)}},\ j=1,2,\cdots, 2^n\}\]
in $H$, where $E_j^{(n)}$ is the partition of $\partial T$ in $2^n$ ``equal" pieces, corresponding to the vertices in the $n$-th level. This partition is invariant with respect to $Aut(T)$, and in particular, with respect to action of $\overline{D}$, $D_{\infty}$ and $\G$. Moreover, $H_n$ naturally embeds into $H_{n+1}$, and \[H=\overline{\bigcup_{n=1}^{\infty}H_n}.\]
Here the bar stands for closure. It is easy to see that $\pi_n$ is unitarily equivalent to the restriction $\pi |_{H_n}$. Similarly, if we denote $\rho_n$ the permutational representation of $D_{\infty}=<a,\ t>$ in $l^2(V_n)$ then $\rho_n$ is unitarily equivalent to $\rho |_{H_n}$.

Let $H_n^{\perp}=H_{n+1}\ominus H_n$ be the orthogonal complement of $H_n$ in $H_{n+1}$ and $\pi_n^{\perp}$ be the restriction $\pi |_{H_n^{\perp}}$. Then
\begin{equation}
\pi={\bf 1}\oplus \bigoplus_{n=0}^{\infty}\pi_n^{\perp},
\end{equation}
where ${\bf 1}$ stands for the trivial representation.
A similar decomposition holds for $\rho$:
\begin{equation}
\rho={\bf 1}\oplus \bigoplus_{n=0}^{\infty}\rho_n^{\perp}.
\end{equation}
In particular, we see that the representations $\pi$ and $\rho$ are sums of finite dimensional representations. The operator recursions (7.5) and (7.6) and the decompositions (7.7) and (7.8) show that the operators $\pi (u)$ and $\rho (t)$ can be presented by the same infinite block matrices as operators on the space
\[H=\C\oplus \bigoplus_{n=0}^{\infty}H_n^{\perp},\]
therefore they coincide. Hence $au$ has infinite order and $\overline{D}$ is isomorphic to $D_{\infty}$. The lemma is proven. \zb

From Lemma 7.3 it follows that the pencils $R_{\rho}(z)=1+z_1\rho(a)+z_2\rho(u)$ and $R_{\pi}(z)=1+z_1\pi(a)+z_2\pi(u)$ coincide. Further, $R_{\rho}(z)$ is invertible in $C^*_{\rho}(D_{\infty})$ if and only if $R_{\pi}(z)$ is invertible in $\C^*_{\rho}(\bar{D})\subset \C^*_{\rho}(\G)$, meaning $P(R_{\rho})=P(R_{\pi})$.
Since $R_{\pi}(z)$ is a linear transform of the pencil $Q(\lb, \mu)=\pi \left(-\lb a+b+c+d-(\mu+1)\right)$, $P(R_{\pi})$ can then be computed through $P(Q)$ using results in \cite{BG}.

Consider the finite dimensional operators
\[Q_n(\lb,\mu)=-\lb a_n+b_n+c_n+d_n-(\mu+1)I_n,\ \ n\geq 0,\]
where \[a_n=\pi_n(a),\ b_n=\pi_n(b),\ c_n=\pi_n(c),\ d_n=\pi_n(d),\]
and $I_n$ is the identity operator of dimension $2^n$. $a_n,\ b_n,\ c_n,\ d_n$ are presented by $2^n\times 2^n$ matrices and satisfiy the following recurrent relations that correspond to the relations in Lemma 7.3:
\begin{align*}
&a_n=\left(\begin{matrix}
                      0 & I_n\\
                      I_n & 0
                     \end{matrix}\right),\ \ \ \ b_n=\left(\begin{matrix}
                                                                      a_{n-1} & 0\\
                                                                      0 & c_{n-1}
                                                                     \end{matrix}\right)\\
&c_n=\left(\begin{matrix}
                      a_{n-1} & 0\\
                      0 & d_{n-1}
                     \end{matrix}\right),\ \ \ \ d_n=\left(\begin{matrix}
                                                                      I_{n-1} & 0\\
                                                                      0 & b_{n-1}
                                                                     \end{matrix}\right),\ \ \ \ n\geq 1.
\end{align*}
As $\pi_n$ is a sub-representation of $\pi$, the joint spectrum for $Q_n$ clearly is contained in the joint spectrum for $Q$. Introduce the polynomials
\[\Phi_0=2-\mu-\lb,\ \ \Phi_1=2-\mu+\lb,\ \ \Phi_2=\mu^2-4-\lb^2,\]
and let
\[\Phi_n=\Phi_{n-1}^2-2(2\lb)^{2^{n-2}},\ \ \ n\geq 3.\]
Set
\[\lb'=\frac{2\lb^2}{4-\mu^2},\ \ \ \mu'=\mu+\frac{\mu \lb^2}{4-\mu^2}\]
and define a map $F:\ \C^2\to \C^2$ by $\lb\to \lb',\ \mu\to \mu'$. Finally, let
\[\Phi_n'=\Phi_n(F(\lb,\mu))=\Phi_n(\lb',\mu').\]
The following proposition is from \cite{BG} Section 4.

\begin{prop}
For all $n\geq 0$, we have

a) $detQ_n(\lb,\mu)=(4-\mu^2)^{2^{n-2}}detQ_{n-1}(F(\lb,\mu))$.

b) $detQ_n=\Phi_0\Phi_1\cdots \Phi_n$.

c) Let $P(Q_n):=\{(\lb,\mu):\ Q_n(\lb,\mu)\ \text{is non-invertible}\}$, then
\[P(Q_n)=\{\Phi_0(\lb,\mu)=0\}\cup \{\Phi_1(\lb,\mu)=0\}\cup \bigcup_{j=1}^{2^{n-1}-1}\{4-\mu^2+\lb^2+4\lb\cos \frac{2\pi j}{2^n}=0\}.\]
\end{prop}
The proof of this proposition is based on the following formula: $\forall n\geq 2,\ \forall k,\ 0\leq k\leq n-2$
\[\Phi_n=\prod_{t=0}^{2^k-1}\big(\Phi_{n-k}-2(2\lb)^2\cos \frac{2\pi (2t+1)}{2^{k+2}}\big).\]
If we set $k=n-2$, then
\begin{align*}
\Phi_n&=\prod_{t=0}^{2^{n-2}-1}\big(\Phi_{2}-4\lb \cos \frac{2\pi (2t+1)}{2^{n}}\big)\\
&=\prod_{t=0}^{2^{n-2}-1}\big(\mu^2-4-\lb^2-4\lb \cos \frac{2\pi (2t+1)}{2^{n}}\big),
\end{align*}
which justifies part c) of the proposition.

Let $H_{x}(\lb,\mu)=4-\mu^2+\lb^2-4\lb x$ and let ${\mathcal Z}(H_x)$ be its zero set in $\C^2$. Then by Proposition 7.4,
 \[\bigcup_{x=\cos \frac{2\pi j}{2^n}} {\mathcal Z}(H_x)\subset P(Q_n)\subset P(Q),\]
where $n\geq 2$ and $1\leq j\leq 2^{n-1}$. Since $P(Q)$ is closed and $\{\cos \frac{2\pi j}{2^n}:\ n\geq 2, 1\leq j\leq 2^{n-1}\}$ is dense in $[-1,\ 1]$ we have
\begin{equation}
\bigcup_{-1\leq x\leq 1} {\mathcal Z}(H_x)\subset P(Q).
\end{equation}
 Converting variables $\lb,\ \mu$ into $z_1=\frac{\lb}{\mu},\ z_2=\frac{-2}{\mu}$ as before, equation $H_{x}(\lb,\mu)=0$ is equivalent to
\[{\mathcal L}_{x}(z)=1-z_1^2-z_2^2-2z_1z_2x=0.\]
Therefore, \[\bigcup_{-1\leq x\leq 1}\{z\in \C^2:\ 1-z_1^2-z_2^2-2z_1z_2x=0\}\subset P(R_{\pi}).\]
Inclusion in the other direction is a consequence of Proposition 2.2, and part b) of the theorem follows.  We note that in the case of real spectrum ($\mu,\ \lb\in {\mathbb R}$) the above fact is proven in \cite{BG} using self-adjointness of $Q$.

Now we prove part a). To this end, we shall need the classification of irreducible unitary representations of $D_{\infty}$ that can be found in \cite{Hal,RS} and \cite{Put}. There are four $1$-dimensional irreducible representations mapping $a,\ t$ of $D_{\infty}$ to the pairs $(1,\ 1),\ (1,\ -1),\ (-1,\ 1)$ or $(-1,\ -1)$, respectively. And, as Example 2.1 indicates, there is a continuum family of inequivalent $2$-dimensional irreducible representations $\rho_{\theta}$. It is not hard to check that $\rho_{\theta}$ is unitarily equivalent to the following representation (also denoted by $\rho_{\theta}$) sending $a$ and $t$ to
\[\left(\begin{matrix}
                      1 & 0\\
                      0 & -1
                     \end{matrix}\right)\ \text{and}\ \left(\begin{matrix}
                                                                      \cos \theta & -\sin \theta\\
                                                                      -\sin \theta & -\cos \theta
                                                                     \end{matrix}\right),\]
respectively, where $0<\theta<\pi$. Observe that if we put $\theta=0$ or $\theta=\pi$ then the corresponding $2$-dimensional representation decomposes as a direct sum of $1$-dimensional representations and in such manner one gets all $1$-dimensional representations. Consider the pencil
\[R_{\theta}(z)=1+z_1\rho_{\theta}(a)+z_2\rho_{\theta}(t)=\left(\begin{matrix}
                      1+z_1+z_2\cos \theta & -z_2\sin \theta\\
                       -z_2\sin\theta & 1-z_1-z_2\cos\theta
                     \end{matrix}\right).\]
As indicated in Example 2.1, $P(R_{\theta})=\{1-z_1^2-z_2^2-2z_1z_2\cos \theta=0\}.$

\begin{prop}
 The Koopman representation $\rho$ of $D_{\infty}$ in $H=L^2(\partial T,\ \mu)$ decomposes as a direct sum
\begin{equation}
\rho={\bf 1}\oplus \mu \oplus \bigoplus_{\theta=\frac{2\pi j}{2^n}}\rho_{\theta},
\end{equation}
where $n\geq 2$ and $1\leq j\leq 2^{n-1}-1$, ${\bf 1}$ is the trivial representation and $\mu$ is the representation given by
\[a\to \left(\begin{matrix}
                      0 & 1\\
                      1 & 0
                     \end{matrix}\right),\ \ \ \ t\to \left(\begin{matrix}
                      1 & 0\\
                      0 & 1
                     \end{matrix}\right).\]
\end{prop}
Note that the $\mu$ above is conjugate to the representation given by $a\to -1,\ t\to 1$. This proposition follows from the equality
\[\bigcup_{n=1}^{\infty}P(R_{\rho_n^{\perp}})=\bigcup_{x=\cos \frac{2\pi j}{2^n}}{\mathcal Z}(H_x),\]
where $R_{\rho_n^{\perp}}(z)=\rho_n^{\perp}(1)+z_1\rho_n^{\perp}(a)+z_2\rho_n^{\perp}(t)$. This was justified in the proof of b). Indeed, we have the decomposition (7.8), where $dim \rho_{n}^{\perp}=2^n$. It allows further decomposition of each $\rho_{n}^{\perp}$ as the the sum of $2$-dimensional irreducible representations. For different values of $x$, where $-1\leq x\leq 1$, the sets
\[\{z:\ 1-z_1^2-z_2^2-2z_1z_2x=0\}\] are different. This implies that each $2$-dimensional irreducible representation of $D_{\infty}$ occurs in (7.10) at most once, and $\rho_{\theta}$ occurs if and only if $\theta=\frac{2\pi j}{2^n}$ for some $n$ and $j$ with $1\leq j\leq 2^{n-1}-1$. This proves the proposition.

To proceed with the proof of Theorem 7.2, we need to say a bit more about weak containment mentioned in Section 2 and introduce the Fell topology on the unitary dual $\widehat{G}$ (the set of equivalence classes of all irreducible unitary representations). If $\pi$ is a unitary representation of a group $G$ then the support of $\pi$, denoted by $supp (\pi)$, is the set of all $\pi'$ in the dual $\widehat{G}$ with $\pi'\prec \pi$. The next proposition is from \cite[Proposition F. 2]{BHV}
\begin{prop}
For a group $G$, any unitary representation $\pi$ of $G$ is weakly equivalent to the direct sum $\oplus \{\pi': \pi'\in supp (\pi)\}$.
\end{prop}
For amenable groups $G$, since the left regular representation $\lb_G$ is maximal with respect to weak containment, the support $supp (\lb_G)=\widehat{G}$.  Applying this fact to the regular representation $\lb_{D_{\infty}}$ of $D_{\infty}$ and using Example 2.1, we get
 \begin{equation*}
\lb_{D_{\infty}} \sim \rho^*:=\bigoplus_{0\leq \theta\leq \pi}\rho_{\theta}.\tag{7.11}
\end{equation*}

For a unitary representation $(\pi, {\mathcal H})$ and an element $x\in {\mathcal H}$, the inner product
\[\phi_x(g)=\langle \pi (g)x,\ x\rangle,\ \ g\in G\]
defines a function on $G$. Functions of this kind are said to be of positive type. \\

{\bf Definition.} For a unitary representation $(\pi, {\mathcal H})$, functions of positive type $\phi_1, \phi_2$, $\cdots, \phi_n$ associated to $\pi$, a compact subset $Q$ of $G$, and $\epsilon  > 0$, let \[W(\pi, \phi_1, \phi_2, \cdots, \phi_n, Q, \epsilon)\] be the set of all unitary representations $\rho\in \widehat{G}$ with the following property:
for each $\phi_i$, there exists a function $\psi$ which is a sum of functions of positive type associated to $\rho$ and such that
\[|\phi_i(x)-\psi(x)|< \epsilon, \ \ \forall x\in Q.\]
The sets $W(\pi, \phi_1, \phi_2, \cdots, \phi_n, Q, \epsilon)$ form a basis for a topology on $\widehat{G}$, called the Fell topology.\\

It is known that if $G$ is separable and locally compact, then $\widehat{G}$ with the Fell topology is separable (\cite{Di}). The following is a consequence of Proposition 7.6 (cf. \cite[Exercise F. 6.4]{BHV}), and it was communicated to us by de la Harpe.
\begin{corr}
If $G$ is a separable locally compact amenable group, and $\{\rho_1, \rho_2, \cdots \}$ is a countable and dense subset in $\widehat{G}$ with respect to the Fell topology, then
\[\lb_G \sim \bigoplus_{i=1}^{\infty} \rho_{i}.\]
\end{corr}

The group $D_{\infty}$ is clearly separable, locally compact and amenable. It is easy to check that the parameterization $\theta\longrightarrow  \rho_{\theta}$ given in Example 2.1 is a homeomorphism from
$(0,\ \pi)$ into $\widehat{D_{\infty}}$ (equipped with the Fell topology) with dense range (\cite{Put}).
Observe that the set of values of the parameter $\theta$ in (7.10) is dense in the interval $[0,\ \pi]$, and hence the set\[\Lambda:=\{\rho_{\theta}:\ \theta=\frac{2\pi j}{2^n}, n\geq 2,\ 1\leq j\leq 2^{n-1}-1\}\]
is dense in $\widehat{D}$. In conclusion, we have by (7.10) and Corollary 7.7 that
\[\lb_{D_{\infty}}\sim \bigoplus_{\rho_{\theta}\in \Lambda}\rho_{\theta} \prec \rho.\]
But since $\lb_{D_{\infty}}$ is maximal, we also have $\rho\prec \lb_{D_{\infty}}$, which concludes that $\lb_{D_{\infty}}$ is weakly equivalent to the Koopman representation $\rho$, and Theorem 7.2 is established. \zb

Theorem 7.1 is an immediate consequence of Theorem 7.2, as $C^*$-algebras generated by weakly equivalent representations of the same group are isomorphic. Further, since $\rho$ is self-similar, $C^*_{\rho}(D_{\infty})$ is self-similar. We end this section with two remarks.\\

{\bf 1}. It is interesting that the element $t\in D_{\infty}$ does not belong to the group $\G$ but belongs to $\C[\G]$. Instead of $\G$ we could use the so-called overgroup $\tilde{\G}=<a,\ \tilde{b},\ \tilde{c},\ \tilde{d}>$ defined by recurrent relations
\[a=\sigma,\ \ \tilde{b}=(a,\ \tilde{c}), \ \ \tilde{c}=(1,\ \tilde{d}),\ \ \tilde{d}=(1,\ \tilde{b}).\]
$\tilde{\G}$ contains $\G$ because $b=\tilde{c}\tilde{d},\ c=\tilde{b}\tilde{d},\ d=\tilde{b}\tilde{c}$. Then one checks that $t=\tilde{b}\tilde{c}\tilde{d}$. The computation of the joint spectrum for the corresponding pencil was provided in \cite{BG} and also can be used in the proof of Theorem 7.2. Observe further that $\tilde{\G}$ is also a group of intermediate growth, but it loses the property to be torsion, as it contains a copy of $D_{\infty}$.

{\bf 2}. The $C^*$-algebra of the infinite cyclic group is also self-similar. This can be seen from the self-similar realization ${\mathbb Z}=<\alpha>$, where
\[\alpha=(1, \alpha)\sigma\]
is an automorphism of a binary rooted tree (called the adding machine or odometer). The spectrum of $\rho(\alpha)$ in the Koopman representation is well-known to be $\{e^{\frac{2\pi ik}{2^n}},\ n\geq 0\}$, and it is dense in the  spectrum of $\lb(\alpha)$ (which is the unit circle ${\mathbb T}$) with respect to the regular representation $\lb$. Hence we have a self-similar realization of $C^*({\mathbb Z})$.
It is an interesting question for which other groups the corresponding group algebras (full or reduced) are self-similar. So far we have only two examples, ${\mathbb Z}$ and $D_{\infty}$.

\section{Dynamics on the joint spectrum}

For the group $\G$ in (0.1), the self-similar relations (7.3) give rise to the map $F:\C^2\to \C^2$ defined in Section 7 as \[(\lb',\mu')=F(\lb,\mu)=\big(\frac{2\lb^2}{4-\mu^2},\ \mu+\frac{\mu \lb^2}{4-\mu^2}\big).\]
$F$ has some interesting properties, and its real-valued dynamics is studied in a series of papers, for example \cite{BG, GNS, GS}. Some of these properties hold true in the complex case.

Let as before $H_{\theta}=4-\mu^2+\lb^2-4\lb\theta$ and $\Psi (\lb,\mu)=\frac{4-\mu^2+\lb^2}{4\lb}$, so we have
\begin{equation}
\Psi(\lb,\mu)=\frac{H_{\theta}(\lb,\mu)+4\lb \theta}{4\lb},
\end{equation}
or equivalently
\begin{equation*}
H_{\theta}(\lb,\mu)=4\lb\Psi(\lb,\mu)-4\lb\theta.\tag{8.1'}
\end{equation*}
\begin{prop}
$F$ is semi-conjugate to the Tchebyshev-von Neumann-Ulam map $\alpha(x)=2x^2-1$.
The semi-conjugation is provided by the map $\Psi$ so that the diagram
\[\begin{tikzcd}
\C^2 \arrow{r}{F} \arrow{d}{\Psi} & \C^2 \arrow{d}{\Psi}\\
\C \arrow{r}{\alpha} & \C
\end{tikzcd}\]
is commutative, i.e. $\Psi(F(\lb,\mu))=2\Psi^2(\lb,\mu)-1$.
\end{prop}
\pf We have by direct computaion that

a)
 \begin{align*}
2\Psi^2(\lb,\mu)-1&=2\frac{(4-\mu^2+\lb^2)^2}{16\lb^2}-1=\frac{(4-\mu^2+\lb^2)^2-8\lb^2}{8\lb^2}\\
&=\frac{16+\mu^4+\lb^4-8\mu^2+8\lb^2-2\mu^2\lb^2-8\lb^2}{8\lb^2}\\
&=\frac{16+\mu^4+\lb^4-8\mu^2-2\mu^2\lb^2}{8\lb^2},
\end{align*}
and

b)
\begin{align*}
\Psi(F(\lb,\mu))&=\frac{4-\mu^2\big(\frac{4-\mu^2+\lb^2}{4-\mu^2}\big)^2+\frac{4\lb^4}{(4-\mu^2)^2}}{\frac{8\lb^2}{4-\mu^2}}\\
&=\frac{4(4-\mu^2)^2-\mu^2(4-\mu^2+\lb^2)^2+4\lb^2}{8\lb^2(4-\mu^2)}\\
&=\frac{4(4-\mu^2)^2-\mu^2(4-\mu^2)^2-2\mu^2(4-\mu^2)\lb^2+\lb^4(4-\mu^2)}{8\lb^2({4-\mu^2})},
\end{align*}
which shows a)=b). \zb

Recall that ${\mathcal Z}(H_{\theta})$ is the zero set of $H_{\theta}$ in $\C^2$, and set $\theta_{1,2}=\pm \sqrt{\frac{1+\theta}{2}}$ (so $\theta_{1,2}$ are preimages of $\theta$ under the map $\alpha(x)=2x^2-1$).
\begin{prop}
$H_{\theta}(F(\lb,\mu))=\frac{1}{4-\mu^2}H_{\theta_1}(\lb,\mu)H_{\theta_2}(\lb,\mu)$. Hence for any $\theta\in \C$ we have the commutative diagram
\[\begin{tikzcd}
F^{-1}({\mathcal Z}(H_{\theta})) \arrow{r}{F} \arrow{d}{\Psi} & {\mathcal Z}(H_{\theta}) \arrow{d}{\Psi}\\
\alpha^{-1}(\theta) \arrow{r}{\alpha} & \theta
\end{tikzcd}\]
\end{prop}
\pf By (8.1'), we have
\begin{align*}
&H_{\theta}(F(\lb,\mu))\\
&=\frac{8\lb^2}{4-\mu^2}(\Psi(F(\lb,\mu))-\theta)\\
&=\frac{8\lb^2}{4-\mu^2}(2\Psi^2(\lb,\mu)-1-\theta)\\
&=\frac{8\lb^2}{4-\mu^2}\big(\frac{(H_{\theta}(\lb,\mu)+4\lb\theta)^2}{8\lb^2}-(1+\theta)\big)\\
&=\frac{1}{4-\mu^2}\big(H_{\theta}(\lb,\mu)+4\lb\theta-2\sqrt{2}\sqrt{1+\theta}\lb\big)\big(H_{\theta}(\lb,\mu)+4\lb\theta+2\sqrt{2}\sqrt{1+\theta}\lb\big)\\
&=\frac{1}{4-\mu^2}(4-\mu^2+\lb^2-2\sqrt{2}\sqrt{1+\theta}\lb)(4-\mu^2+\lb^2+2\sqrt{2}\sqrt{1+\theta}\lb)\\
&=\frac{1}{4-\mu^2}H_{\theta_1}(\lb,\mu)H_{\theta_2}(\lb,\mu).
\end{align*}
\zb

$F$'s fixed points can be easily determined. Consider the equations
\[\frac{2\lb^2}{4-\mu^2}=\lb,\ \ \ \mu+\frac{\mu\lb^2}{4-\mu^2}=\mu.\]
Clearly, every point in the plane $\{\lb=0\}$ is a fixed point. If $\lb\neq 0$, then $\mu=0$ by the second equation, which in turn implies $\lb=2$. So the set of fixed points is
\[S=\{\lb=0\}\cup\{(2,0)\}.\]
One sees that $(2,0)$ is a zero of $\Phi_0$ (Section 7), and hence it is in $P(Q)$ by Proposition 7.4. But the plane $\{\lb=0\}$ is not in $P(Q)$.

Parallel analysis can be done for $D_{\infty}$. Using the recursion relations in (7.2), we can write
\[R(z)=z_0+z_1a+z_2t=\begin{pmatrix}
                                      z_0+z_2a & z_1\\
                                      z_1 & z_0+z_2t
                                      \end{pmatrix}.\]
Assume $z_0^2\neq z_2^2$, then $z_0+z_2a$ is invertible and its inverse is $(z_0-z_2a)(z_0^2-z_2^2)^{-1}$.
Hence by (1.4) $R(z)$ is invertible if and only if $z_0+z_2t-z_1^2(z_0-z_2a)(z_0^2-z_2^2)^{-1}$ is invertible, or equivalently if and only if
\begin{equation}
\frac{z_0(z_0^2-z_1^2-z_2^2)}{z_0^2-z_2^2}+\frac{z_1^2z_2}{z_0^2-z_2^2}a+z_2t
\end{equation}
is invertible. We define $F_1:\C^3\to \C^3$ by
\begin{equation}
F_1(z_0,z_1,z_2)=\big(z_0(z_0^2-z_1^2-z_2^2),z_1^2z_2,(z_0^2-z_2^2)z_2\big):=(z_0',z_1',z_2').
\end{equation}
Since $P(R)$ is symmetric in $z_1$ and $z_2$, we can also define
\[F_2(z_0,z_1,z_2)=\big(z_0(z_0^2-z_1^2-z_2^2),z_1z_2^2,(z_0^2-z_1^2)z_1\big).\]
But since the properties are parallel, we shall only focus on $F_1$.
\begin{thm}
The joint spectrum $P(R)$ in Theorem 1.1 is invariant under $F_1$.
\end{thm}
\pf We need to show that $F_1(P(R))\subset P(R)$. First, with the assumption $z_1z_2\neq 0$ and $z_0^2\neq z_2^2$, we compute $\frac{(z_0')^2-(z_1')^2-(z_2')^2}{2z_1'z_2'}$. One checks that
\begin{align*}
&(z_0')^2-(z_1')^2\\
&=(z_0'-z_1')(z_0'+z_1')\\
&=\big(z_0(z_0^2-z_1^2-z_2^2)-z_1^2z_2\big)\big(z_0(z_0^2-z_1^2-z_2^2)+z_1^2z_2\big)\\
&=(z_0+z_2)\big(z_0(z_0-z_2)-z_1^2\big)(z_0-z_2)\big(z_0(z_0+z_2)-z_1^2\big),
\end{align*}
and therefore substituting $z_1'$ by $z_1^2z_2$ and $z_2'$ by $(z_0^2-z_2^2)z_2$, we have
\begin{align*}
&\frac{(z_0')^2-(z_1')^2-(z_2')^2}{2z_1'z_2'}\\
&=\frac{\big(z_0(z_0-z_2)-z_1^2\big)\big(z_0(z_0+z_2)-z_1^2\big)-(z_0^2-z_2^2)z_2^2}{2z_1^2z_2^2}\\
&=\frac{z_0^4+z_1^4+z_2^4-2z_0^2z_1^2-2z_0^2z_2^2}{2z_1^2z_2^2}\\
&=\frac{(z_0^2-z_1^2-z_2^2)^2-2z_1^2z_2^2}{2z_1^2z_2^2}\\
&=2\big(\frac{z_0^2-z_1^2-z_2^2}{2z_1z_2}\big)^2-1.
\end{align*}
Now we recall the Tchebyshev-von Neumann-Ulam map $\alpha(x)=2x^2-1,\ \ x\in \C$. It is easy to check that
$\alpha$ maps $[-1,\ 1] \to [-1,\ 1]$ and $\C\setminus [-1,\ 1]\to \C\setminus [-1,\ 1]$. Then it follows from Theorem 1.1 that $F_1$ maps $P(R)\setminus \left(\{z_1z_2=0\}\cup \{z_0^2=z_2^2\}\right)$ into $P(R)$.

Now we look at the case $z_1z_2=0$. By symmetry, we only deal with the case $z_1=0$. If
$z\in P(R)\cap \{z_1=0\}$, then by Theorem 1.1, we have $z_0^2-z_2^2=0$, and hence $F_1(z)=(0, 0, 0)\in P(R)$.

Finally, if  $z_0^2-z_2^2=0$, then $z_0= \pm z_2$. Since in this case
 \[z'=F_1(z)=z_1^2\left(-z_0,\ z_2,\  0\right),\]
we have \[R(z')=z_0'I+z_1'\lb_{D_{\infty}}(a)+z_2'\lb_{D_{\infty}}(t)=z_1^2(-z_0I+z_2\lb_{D_{\infty}}(a)).\]
Since the classical spectrum $\sigma(\lb_{D_{\infty}}(a))=\{\pm 1\}$, we see that $R(z')$ is not invertible. This completes the proof. \zb

Now consider \[{\mathcal L}_{\theta}=z_0^2-z_1^2-z_2^2-2z_1z_2\theta,\]
and for $\theta \in \C$, we denote $\alpha(\theta)$ by $\theta'$. Then from above computations we have
\[\frac{(z_0')^2-(z_1')^2-(z_2')^2}{2z_1'z_2'}-\theta'=2\big(\frac{z_0^2-z_1^2-z_2^2}{2z_1z_2}\big)^2-2\theta^2\]
which implies
\begin{align}
{\mathcal L}_{\theta'}(z')&=4z_1'z_2'\frac{{\mathcal L}_{\theta}(z){\mathcal L}_{-\theta}(z)}{4z_1^2z_2^2}\\
&=(z_0^2-z_2^2){\mathcal L}_{\theta}(z){\mathcal L}_{-\theta}(z).
\end{align}
Since ${\mathcal L}_{\theta}(z)$ and $F_1(z)$ are holomorphic everywhere, and (8.5) hold on the open set
$\{ z\in \C^3:\ z_1z_2\neq 0,\ z_0^2\neq z_2^2\}$, it holds on the entire $\C^3$. So parallel to Proposition 8.2 we have
\begin{prop}
${\mathcal L}_{\theta}(F_1(z))=(z_0^2-z_2^2){\mathcal L}_{\theta_1}(z){\mathcal L}_{\theta_2}(z), \ \ z\in \C^3.$
\end{prop}

Theorem 8.3 indicates that $P(R)$ is an invariant set for $F_1$. It is interesting to observe that the three hyperplanes $S_0=\{z_0=0\},\ S_1=\{z_1=0\}$ and $S_2=\{z_2=0\}$ are also invariant for $F_1$. $F_1$'s fixed points can also be determined. Consider the equations
\begin{equation}
z_0(z_0^2-z_1^2-z_2^2)=z_0,\ \ z_1^2z_2=z_1,\ \ (z_0^2-z_2^2)z_2=z_2.
\end{equation}

Here are two cases.

1. If $z_0\neq 0$, then $z_0^2-z_1^2-z_2^2=1$, and hence the third equation becomes $(1+z_1^2)z_2=z_2$, which, by the second equation, implies $z_1=0$. It is then easy to see that every point in the set
\[S_3=\{z_1=0,\ \ z_0^2-z_2^2=1\}\]
is a fixed point for $F_1$.

2. If $z_0=0$, then the third equation becomes $-z_2^3=z_2$. If $z_1\neq 0$, then $z_2\neq 0$ by the second equation, and hence $z_2=\pm i$ and $z_1=1/z_2$. If $z_1=0$, then by the third equation $z_2=0,\ \pm i$. Therefore every point in
\[S_4=\{(0,i,-i),(0,-i,i),(0,0,0),(0,0,\pm i)\}\]
is a fixed point for $F_1$. Checking with Theorem 1.1, one sees that
\[S_3\cup \{(0,0,\pm i)\}\subset P^c(R)\] and \[\{(0,i,-i),(0,-i,i),(0,0,0)\}\subset P(R).\]

In addition, $F_1$'s Jacobian can be computed as
\begin{align*}
J(F_1)(z)&=det \begin{pmatrix}
                     3z_0^2-z_1^2-z_2^2 & 0 & 2z_0z_2\\
                         -2z_0z_1 & 2z_1z_2  & 0\\
                         -2z_0z_2 & z_1^2 & z_0^2-3z_2^2
                                      \end{pmatrix}\\
&=(3z_0^2-z_1^2-z_2^2)2z_1z_2 (z_0^2-3z_2^2)+2z_0z_2( -2z_0z_1^3+4z_0z_1z_2^2)\\
&=2z_2\big(3z_0^4z_1-6z_0^2z_1z_2^2+3z_1^3z_2^2+3z_1z_2^4-3z_0^2z_1^3\big)\\
&=6z_1z_2(z_0^2-z_2^2)(z_0^2-z_1^2-z_2^2),
\end{align*}
from which one easily determines $F_1$'s singularities. Similar to the studies in \cite{BG, GNS, GS}, there is much more one can say about the dynamics of $F$ and $F_1$.

\section{Concluding remarks}

The objective of this paper is to do an exploration on the idea of joint spectrum in group theory. It shows that in the case of dihedral group the projective joint spectrum $P(R)$ can be computed, it reflects the structure of $D_{\infty}$, it is a good invariant for the representations, and it is a measurement for weak containment. These results nurture an anticipation that joint spectrum may have an interesting role to play in group theory.

\vspace{5mm}


\begin{thebibliography}{10}

\bibitem{AJ} M. Andersson and J. S\"{o}strand, {\em Functional calculus for non-commuting operators with real spectra via an iterated cauchy formula}, J. Funct. Anal. 210 (2004), No.2, 341-375.

\bibitem{At} F. V. Atkinson, {\em Multiparameter eigenvalue problems}, Academic Press, New York and London, 1972.

\bibitem{Av} W. Arveson, {\em An invitation to $C^*$-algebra}, Springer-Verlag, 1981.

\bibitem{BCY} J. Bannon, P. Cade and R. Yang, {\em On the spectrum of operator-valued entire functions}. Illinois J. of Mathematics 55 No.4 (2011).

\bibitem{BG} L. Bartholdi and R. Grigorchuk, {\em On the spectrum of Hecke type operators related to some fractal groups}, Tr. Mat. Inst. Steklova 231 (2000), Din. Sist., Avtom. i Beskon. Gruppy, 5--45; translation in Proc. Steklov Inst. Math. 2000, no. 4 (231), 1–41.

\bibitem{BHV} B. Bekka, de la Harpe and A. Valette, {\em Kazhdan's property (T)}, New Mathematical Monographs 11, Cambridge University Press, Cambridge, 2008.

\bibitem{BG..S} I. Bondarenko, R.  Grigorchuk, R. Kravchenko, Y. Muntyan, V. Nekrashevych, D. Savchuk and Z. \v{S}uni\'{c}, {\em Groups generated by 3-state automata over a 2-letter alphabet}, I. S\~{a}o Paulo J. Math. Sci. 1 (2007), no. 1, 1–39.

\bibitem{CY} P. Cade and R. Yang, {\em Projective spectrum and cyclic cohomology}. J. of Funct. Analy. Vol. 265 No. 9 (2013).

\bibitem{CSZ} I. Chagouel, M. Stessin and K. Zhu, {\em Geometric spectral theory for compact operators}, Trans. Amer. Soc. 368 (2016), No. 3, 1559-1582.

\bibitem{Da} K. R. Davidson, {\em $C^*$-algebras by example}, Fields Institute Monographs,
A.M.S, 1996.

\bibitem{De} C. Deninger, {\em Mahler measure and Fuglede-Kadison determinant}, M\"{u}ster J. Math (2) (2009), 45-63.

\bibitem{Di} J. Dixmier, {\em Les $C^*$-alg\`{e}bres et leurs repr\'{e}sentations}, Gauthier-Villars, 1969.

\bibitem{DY} R. G. Douglas and R. Yang, {\em Hermitian geometry on resolvent sets}, preprint.

\bibitem{DG} A. Dudko and R. Grigorchuk, {\em On spectra of Koopman, groupoid and quasi-regular representations}, preprint.

\bibitem{EW} G. Everest and T. Ward, {\em Height of polynomials and entrophy in algebraic dynamics}, Universitext, Springer-Verlag London Limited, 1999.

\bibitem{Fa} A. Fainshtein, {\em Taylor joint spectrum for family of operators generating nil-potent Lie algebras}, J. Oper. Theory 29 (1993), 3-27.

\bibitem{FK} B. Fuglede and R. Kadison, {\em Determinant theory in finite factors}, Annals of Math. 55 (1952), 520-530.

\bibitem{Gr80} R. Grigorchuk, {\em  On Burnside's problem on periodic groups}, (Russian) Funktsional. Anal. i Prilozhen. 14 (1980), no. 1, 53–54.

\bibitem{Gr84} R. Grigorchuk, {\em Construction of p-groups of intermediate growth that have a continuum of factor-groups}, (Russian) Algebra i Logika 23 (1984), no. 4, 383–394, 478.

\bibitem{Gr89} R. Grigorchuk, {\em On the Hilbert-Poincare series of graded algebras that are associated with groups}, Mat. Sb. 180 (1989), No. 2, 207-225; translation in Math. USSR-Sb. 66 (1990), No. 1, 211-229.

\bibitem{Gr05} R. Grigorchuk, {\em Solved and unsolved problems around one group}, Infinite groups: geometric, combinatorial and dynamical aspects, 117-218, Progr. Math., 248, Birkhauser, Basel, 2005.

\bibitem{GN07} R. Grigorchuk and V. Nekrashevych, {\em Self-similar groups, operator algebras and Schur complement}, J. Mod. Dyn. 1 (2007), no. 3, 323-370.

\bibitem{GNS} R. Grigorchuk, V. Nekrashevich and V. Sushchanski\v{i}, {\em  Automata, dynamical systems, and groups}, (Russian) Tr. Mat. Inst. Steklova 231 (2000), Din. Sist., Avtom. i Beskon. Gruppy, 134--214; translation in Proc. Steklov Inst. Math. 2000, no. 4 (231), 128–203.

\bibitem{GS} R. Grigorchuk and Z. \v{S}uni\'{c}, {\em Schreier spectrum of the Hanoi Towers group on three pegs}, Proc. of Symposia in Pure Math. Vol. 77, 2008.

\bibitem{Hal} P. R. Halmos, {\em Two subspaces}, Trans. Amer. Math. Soc. 144 (1969), 381-389.

\bibitem{Ha} P. de la Harpe, {\em The Fuglede-Kadison determinant, theme and variations}, Proc. of the National Academy of Sci. of U.S.A, Vol 110, no. 40 (2013), 15864-15877.

\bibitem{Ha72} R. Harte, {\em Spectral mapping theorems}, Proc. Royal Irish Acad., 72 A (1972), 89-107.

\bibitem{Ha73} R. Harte, {\em Spectral mapping theorems for quasi-commuting systems}, Proc. Royal Irish Acad., 73 A (1973), 7-18.

\bibitem{HS} P. de la Harpe and G. Skandalis, {\em D\'{e}terminant associ\'{e} \`{a} une trace sur une alg\'{e}bre de Banach} (French), Ann. Inst. Fourier (Grenoble) 34 (1984), No. 1, 241-260.

\bibitem{HY} W. He and R. Yang, {\em Projective spectrum and kernel bundle}. Sci. China. Math. Vol. 57 (2014), 1-10.

\bibitem{Ho} L. H\"{o}rmander, {\em An introduction to complex analysis in several variables}, 3rd ed., North Holland, Amsterdam, 1990.

\bibitem{Je} B. Jefferies, {\em Spectral properties of noncommuting operators}, Lecture Notes in Math, 1843, Springer-Verlag, Berlin 2004.

\bibitem{Li} H. Li, {\em Compact group automorphisms, addition formulas and Fuglede-Kadison determinants}, Ann. of Math. (2) 176 (2012), no. 1, 303-347.

\bibitem{LS} T. Lu and S. Shiou, {\em Inverse of $2\times 2$ block matrices}, Comp. and Math. with Appl. 43 (2002), 119-129.

\bibitem{Ne05} V. Nekrashevych, {\em Self-similar groups}, Mathematical Survey and Monographs, A.M.S Providence, RI, 2005.

\bibitem{Ne} V. Nekrashevych, {\em Periodic groups from minimal actions of the infinite dihedral group}, arXiv:1601.01033 (Jan. 2016).

\bibitem{Pe} G. K. Pedersen, {\em Measure theory in $C^*$-algebras II}, Math. Scand. 22 (1968), 63-74.

\bibitem{Put} I. Putnam, {\em Lecture Notes on $C^*$-algebras}, preprint.

\bibitem{RS} I. Raeburn and A. Sinclair, {\em The $C^*$-algebras generated by two projections}, Math. Scand. 65 (1989), 278-290.

\bibitem{Ru} W. Rudin, {\em Real and complex analysis}, McGraw-Hill Series in Higher Mathematics, McGraw-Hill, 1996.

\bibitem{Sc} K. Schmidt, {\em Dynamical systems of algebraic origin}, Modern Birkhauser Classics, Birkhauser/Springer Basel AG, Basel, 1995.

\bibitem{Sl} B. D. Sleeman, {\em Multiparameter spectral theory in Hilbert space}, Research Notes in Mathematics, Vol. 22, Pitman, London 1978.

\bibitem{SYZ} M. Stessin, K. Zhu and R. Yang, {\em Analyticity of a joint spectrum and a multivariable analytic Fredholm theorem}. New York J. Math. 17A (2011), 39-44.

\bibitem{Ta70} J. L. Taylor, {\em A joint spectrum for several commutative operators}, J. Functional Analysis 6 (1970) 172--191.

\bibitem{Ta72} J. L. Taylor, {\em A general framework for multi-operator functional calculus}, Advances Math. 9 (1972), 137--182.

\bibitem{Vi} V. Vinnikov, {\em Determinantal representation of algebraic curves}, Linear algebra in signals, systems and control (Boston, MA, 1986), 73-99, SIAM, Philadelphia, PA, 1988.

\bibitem{Ya} R. Yang, {\em Projective spectrum in Banach algebras}, J. Topol. and Analy. 1 (2009), No. 3, 289-306.

\bibitem{ZKKP} M. Zaidenberg, S. Krein, P. Kuchment and A. Pankov, {\em Banach bundles and linear operators}, Russian Math. Surveys, 30:5 (1975), 115-175.

\end{thebibliography}
\end{document}